\documentclass[a4paper,10pt]{article}
\usepackage[utf8]{inputenc}
\usepackage[T1]{fontenc}
\usepackage[utf8]{inputenc}
\usepackage{amsthm}
\usepackage{amsmath,amsfonts,mathrsfs}
\usepackage{hyperref}
\usepackage{xspace}
\usepackage{graphicx}
\usepackage[rmargin=3cm,lmargin=3cm,tmargin=2cm,bmargin=3cm]{geometry}
\usepackage{amssymb}
\usepackage{textcomp}
\usepackage{wasysym}
\usepackage{xcolor}
\newtheorem{thm}{Theorem}[section]

\newtheorem{lem}[thm]{Lemma}
\newtheorem{prop}{Proposition}[section]
\theoremstyle{definition}

\newtheorem{exam}{Example}[section]
\newtheorem{rem}{Remark}[section]

\numberwithin{equation}{section}

\newcommand{\al}{\alpha}
\newcommand{\si}{\sigma}

\newcommand{\ga}{\gamma}

\newcommand{\be}{\beta}

\newcommand{\de}{\delta}
\newcommand{\De}{\Delta}

\DeclareMathOperator{\tr}{tr}
\DeclareMathOperator{\di}{div}

\DeclareMathOperator{\grad}{grad}
\newcommand{\Rom}{R}

\newcommand{\R}{\mathbb{R}}
\newcommand{\V}{\mathcal{V}}

\newcommand{\Hi}{\mathcal{H}}

\newcommand{\Real}{\mathbb{R}}
\newcommand{\Vol}{Vol}

\newcommand{\X}{\mathfrak{X}}

\title{Harmonic vector fields on extended 3-dimensional Riemannian Lie groups}
\author{F. Koudjo\footnote{Institut de Mathematiques et de Sciences Physiques, Porto-Novo, Benin,\ ferdinand.koudjo@imsp-uac.org}, 
E. Loubeau\footnote{Univ. Brest, CNRS UMR 6205, LMBA, F-29238 Brest, France,\ Eric.Loubeau@univ-brest.fr} 
and  L. Todjihounde\footnote{Institut de Mathematiques et de Sciences Physiques, Porto-Novo, Benin,\ leonardt@imsp-uac.org}}

\date{}

\begin{document}
	
	\maketitle
	
	\begin{abstract} Given two Riemannian manifolds $(B,g_B)$ and $(F,g_F)$, we  give  harmonicity conditions for  vector fields on the Riemannian warped product $B\times_fF$, with $f:B \longrightarrow ]0,+\infty[$, using a characteristic variational condition. Then, we apply this to the case $B=\mathbb{R}$ and $F$ is a three-dimensional connected Riemannian Lie group $G$ equipped with a left-invariant metric, to determine harmonic vector fields on $\R\times_fG$. We give examples of harmonic vector fields on $G$ which are not left-invariant and determine harmonic vector fields on $\R\times_fG$. We conclude with some examples of vector fields on $\R\times_fG$ which are harmonic maps. 
\end{abstract}
	

	\section{Introduction}
	
	One of the most studied objects in Differential Geometry is the energy functional of a map\\ $\varphi:(M^m,g)\rightarrow (N^n, h)$ between
	Riemannian manifolds of dimensions $m$ and $n$, respectively, given by
	\begin{equation*}
	E(\varphi)=\int_{D}{e(\varphi)v_{g}}.
	\end{equation*}
	where $D$ is a compact domain of $M$, $e(\varphi):M \rightarrow [0,\infty[$ the energy density of $\varphi$ defined by \begin{equation*}
	e(\varphi)(x)=\dfrac{1}{2}\|d\varphi_{x}\|^{2}=\dfrac{1}{2}\sum_{i=1}^{m}{h(d\varphi_{x}(e_{i}), d\varphi_{x}(e_{i}))},
	\end{equation*}
	for $x\in M,\{e_{i}\}_{i=1}^{m}$ an orthonormal basis of$T_{x}M$ and  $d\varphi_{x}$ the differential of the map $\varphi$ at the	point $x$ (\cite{Baird03,EellsSampson64}).
	
	Denote by $C^{\infty}(M, N)$ the space of  smooth maps from $M$ to $N$, $\nabla^{\varphi}$ the connection
	of the vector bundle $\varphi^{-1}TN$ induced from the Levi-Civita connection $\bar{\nabla}$ of $(N, h)$ and $\nabla$ the Levi-Civita connection of $(M, g)$.
	
	A map $\varphi: (M, g)\rightarrow (N, h)$ is said to be \textit{harmonic} if it is a critical point of the energy functional
	$E(.;D):C^{\infty}(M, N)\rightarrow \Real$ for any compact domain $D$. It is well-known (\cite{EellsSampson64}) that the map $\varphi:
	(M, g)\rightarrow (N, h)$ is harmonic if and only if
	\begin{equation}
	\tau(\varphi)=\tr(\nabla
	d\varphi)=\sum_{i=1}^{m}\big\{\nabla_{e_{i}}^{\varphi}d\varphi(e_{i})-d\varphi(\nabla_{e_{i}}{e_{i}})\big\}=0.\end{equation}
	The	equation $\tau(\varphi)=0$ is called the \emph{harmonic equation}.	Denote by $(TM, g_{S})$ the tangent bundle of $(M, g)$ equipped with the Sasaki metric $g_{S}$ (cf. section 2). A vector field $X$ on $M$ determines a map from $(M, g)$ into $(TM, g_{S})$ (i.e. a section of $TM$), embedding $M$ into its tangent bundle.
	Its \emph{energy} $E(V)$  is given by
	\begin{equation*}
	E(V)=\dfrac{m}{2}\Vol(D)+\frac{1}{2}\int_{D}{\|\nabla V\|^{2}v_{g}}\approx Vol(D)+E^v(V).
	\end{equation*}
	
	It was shown in \cite{Ishiha79} and \cite{Nouhaud77} that if $M$ is compact and a vector field $X$ is a harmonic map from $(M, g)$ into $(TM, g_{S})$, then $X$ must be parallel. This rigidity may be overcome for vector fields of unit length. A unit vector field $V$ on $(M, g)$ is a section of the unit tangent sphere bundle $T_{1}M,$ i.e. $V:(M, g)\rightarrow (T_{1}M, g_{S})$, where $T_{1}M$ carries the metric induced from $g_{S}$ \cite{Wiegm95}. Critical points of $E_{1}$, with respect to variations through  vector fields (respectively unit vector fields) are called \emph{harmonic  vector fields} (\emph{resp. harmonic unit vector fields}).
	The corresponding critical point conditions have been determined in \cite{Wiegm95} and \cite{Wood97}. It should be pointed out that a harmonic  vector field determines a harmonic map when an additional condition involving the curvature is satisfied (\cite{GilMedr01,HanYim98}).
		
	\section{Preliminaries}
	
	\subsection{The tangent bundle and the unit tangent sphere
		bundle}
	
	Let $(M, g)$ be an $m$-dimensional Riemannian manifold and $\nabla$
	the associated Levi-Civita connection. Its Riemann curvature tensor
	$R$ is defined by
	\begin{equation*}
	R(X, Y)Z =\nabla_{X}\nabla_{Y}Z - \nabla_{Y}\nabla_{X}Z -
	\nabla_{[X, Y]}Z
	\end{equation*}
	for all vector fields $X, Y$ and $Z$ on $M$. The tangent bundle of $(M, g)$, denoted by $TM$, consists of
	pairs $(x, u)$ where $x$ is a point in $M$ and $u$ a tangent vector
	to $M$ at $x$. The mapping $\pi: TM \rightarrow M: (x, u)\mapsto x$ is
	the natural projection from $TM$ onto $M$. The tangent space $T_{(x, u)}TM$ at a point $(x,
	u)$ in $TM$ is the direct sum of the vertical subspace $\V_{(x,
		u)}=Ker(d\pi|_{(x, u)})$ and the horizontal subspace $\Hi_{(x,
		u)}$, with respect to the Levi-Civita connection $\nabla$ of $(M, g)$:
	\begin{equation*}
	T_{(x, u)}TM=\Hi_{(x, u)}\oplus \V_{(x, u)}.
	\end{equation*}
	For any vector $w\in T_{x}M$, there exists a unique vector $w^{h}\in
	\Hi_{(x, u)}$ at the point $(x, u)\in TM$, called the
	\emph{horizontal lift} of $w$ to $(x, u)$, such that $d\pi(w^{h})=w$
	and a unique vector $w^{v}\in \V_{(x, u)}$, the
	\emph{vertical lift} of $w$ to $(x, u)$, such that $w^{v}(df)=w(f)$
	for all functions $f$ on $M$. Hence, every tangent vector
	$\bar{w}\in T_{(x, u)}TM$ can be decomposed as
	$\bar{w}=w_{1}^{h}+w_{2}^{v}$ for uniquely determined vectors
	$w_{1}, w_{2} \in T_{x}M$. The \emph{horizontal} (resp.
	\emph{vertical}) \emph{lift} of a vector field $X$ on $M$ to $TM$ is
	the vector field $X^{h}$ (resp. $X^{v}$) on $TM$ whose value
	at the point $(x, u)$ is the horizontal (respectively, vertical)
	lift of $X_{x}$ to $(x, u)$.
	
	The tangent bundle $TM$ of a Riemannian manifold $(M, g)$ can be
	endowed in a natural way with a Riemannian metric $g_{S}$, the
	\emph{Sasaki metric}, depending only on the Riemannian structure $g$
	of the base manifold $M$. It is uniquely determined by
	\begin{equation}\label{Eq: Sasaki metric}
	\begin{array}{ll}
	g_{S}(X^{h}, Y^{h})=g_{S}(X^{v}, Y^{v})=g(X, Y)\circ \pi, &g_{S}(X^{h}, Y^{v})=0\\
	\end{array}
	\end{equation}
	for all vector fields $X$ and $Y$ on $M$. More intuitively, the
	metric $g_{S}$ is constructed in such a way that the vertical and
	horizontal subbundles are orthogonal and the bundle map $\pi: (TM,
	g_{S})\mapsto (M, g)$ is a Riemannian submersion.
	We denote by $\X(M)$ the set of globally defined vector fields on the base manifold $(M, g)$. In the sequel, we concentrate on the map $V: (M, g)\rightarrow (TM,g_{S})$(resp. $V: (M, g)\rightarrow (T_{1}M,
	g_{S})$). The tension field $\tau_{1}(V)$ of $V: (M, g)\rightarrow (TM,g_{S})$ is given by \cite{GilMedr01}
	\begin{equation}\label{Eq: tension field of V}
	\tau(V)=(-S(V))^{h}+(-\bar{\Delta}V)^{v},
	\end{equation} 
	where $\{e_{i}\}_{i=1}^{m}$ is a local orthonormal frame field of $(M, g)$, $\displaystyle S(V)=\sum_{i=1}^{m}{\Rom(\nabla_{e_{i}}{V}, V)e_{i}}$ and \\$\displaystyle \bar{\Delta}V=\nabla^*\nabla V=\sum_{i=1}^{m}\{\nabla_{\nabla_{e_{i}}{e_{i}}}V-\nabla_{e_{i}}\nabla_{e_{i}}V\}$. Consequently, $V$ defines a harmonic map from  $ (M,g)$ to $ (TM,g_s)$  if and only if 
	\begin{eqnarray}
	\tr[R(\nabla_.V,V).]=0 \quad\text{ and} \quad\nabla^*\nabla V=0.
	\end{eqnarray}
	
	A smooth vector field $V$ is said to be a harmonic section if and only if it is a critical point of $E^v$ where $ E^v$ is the vertical energy. The corresponding Euler-Lagrange equation is given by $$\nabla^*\nabla V=0.$$

	\subsection{ Warped Products}
	
	Let $(B^m,g_B)$ and $(F^n,g_F)$ be Riemannian manifolds with $f: B\longrightarrow \text{ }]0,+\infty[$ a smooth function on B. The warped product $M =B\times_fF$ is the product
	manifold $B\times F$ equipped with the metric
	$g = \pi^*(g_B)\oplus(f\circ\pi)^2\sigma^*(g_F )$;
	where $\pi:M\longrightarrow B$ and $ \sigma: M\longrightarrow F$ are the usual projections. Then
	$(B, g_B)$ is called the base, $(F, g_F)$ is the fiber and $f$ the warping function of the warped product, $\pi^{-1}(p)=\{p\}\times F$ are the fibers and $\si^{-1}(q)=B\times\{q\}$ the leaves. The vectors tangent to leaves are called horizontal and those tangent to the fibers vertical, hence \begin{eqnarray*}
		T_{(p,q)}(B\times F)&=&T_{(p,q)}(\{p\}\times F)\oplus T_{(p,q)}(B\times \{q\})\\&=&T_{(p,q)}(\{p\}\times F)\oplus T_{(p,q)}(\{p\}\times F)^{\perp}\\&=&T_{(p,q)}(B\times \{q\})^\perp\oplus T_{(p,q)}(B\times \{q\})\end{eqnarray*}for all $(p,q)\in M$.\\
	Vector fields $X$ on $B\times F$ are horizontal if $$d\pi_{(p,q)}(X)=X_p\quad\text{and }\quad d\si_{(p,q)}(X)=0$$
	that is  $$X\in T_{(p,q)}(B\times \{q\})^\perp \quad\text{and }\quad d\pi_{(p,q)}(X)=X_p.$$
	
	Vector fields $X$ on $B\times F$ are vertical  if $$d\pi_{(p,q)}(X)=0\quad\text{and }\quad d\si_{(p,q)}(X)=X_q$$
	that  is $$X\in T_{(p,q)}(\{p\}\times F)^\perp \quad\text{and }\quad d\si_{(p,q)}(X)=X_p.$$ 
	
	If $X\in T_pB$ and $q \in F$ then the horizontal  lift of $X$ to  $(p,q)$ is the unique vector $X^*$ in $T_{(p,q)}(B\times F)$ such that
	$d\pi_{(p,q)}(X^*)=X_p \quad\text{and }\quad d\si_{(p,q)}(X)=0$
	
	If $X\in T_p(F)$ and $q \in F$ then the vertical lift of $X$ to  $(p,q)$ is the unique vector $X^*$ in $T_{(p,q)}(B\times F)$ such that
	$d\pi_{(p,q)}(X)=0 \quad\text{and }\quad d\si_{(p,q)}(X)=X_q$
	\begin{lem}[\cite{Oneil01}]
		Let $(B^m,g_B)$ and $(F^n,g_F)$ be Riemannian manifolds with $f: B\longrightarrow \text{ }]0,+\infty[$ a smooth function on B. Let $X_1,Y_1$ be vector fields on $B$ and $X_2,Y_2$ vector fields on $F$. Let  $\nabla,\nabla^B,\nabla^F$ be the Levi-Civita connections of $(M,g)$, $(B,g_B)$ and $(F,g_F)$ respectively, then \begin{enumerate}
			\item $\grad^M(f\circ\pi)=(\grad^B(f),0)$;
			\item $\nabla_{(X_1,0)}(Y_1,0)=(\nabla^B_{X_1}Y_1,0)$;
			\item $\nabla_{(0,X_2)}(0,Y_2)=(0,\nabla^F_{X_2}Y_2)-fg^F(X_2,Y_2) \grad^B(f\circ\pi)\\ \\\text{ }\quad\text{ }\quad\text{ }\quad\qquad=(0,\nabla^F_{X_2}Y_2)-\dfrac{1}{f}g((0,X_2),(0,Y_2))\grad^B(f\circ\pi)$
			\item $\nabla_{(0,X_2)}(Y_1,0)=\dfrac{Y_1(f)}{f}(0,X_2)$
			\item $\nabla_{(X_1,0)}(0,Y_2)=\dfrac{X_1(f)}{f}(0,Y_2)$;
			\item  $\grad^M(h\circ\si)=\dfrac{1}{f^2}(0,\grad^Fh)$, \quad for \quad $h:F\longrightarrow \mathbb{R}$.
		\end{enumerate}

	\end{lem}

	\begin{lem}[\cite{Oneil01}]
		Let $(B^m, g_B)$ and $(F^n, g_F)$ be Riemannian manifolds, with $f: B\longrightarrow ]0,+\infty[$ a smooth function on B, $X_1,Y_1,Z_1$ vector fields on $B$, $X_2,Y_2,Z_2$ vector fields on $F$, $\nabla,\nabla^B$ the Levi Civita connections of $(M,g)$, $(B,g_B)$ respectively, and	$R,R^B,R^F$ the Riemmanian curvature tensors on  $(M,g)$, $(B,g_B)$ and $(F,g_F)$ respectively,	then \begin{enumerate}
			\item $R\left((X_1,0),(Y_1,0)\right)(Z_1,0)=\left(R^B(X_1,Y_1)Z_1,0\right)$\\
			\item $R\left((X_1,0),(Y_1,0)\right)(0,Z_2)=0$\\
			\item  $R\left((0,X_2),(0,Y_2)\right)(Z_1,0)=0$\\
			\item$R\left((0,X_2),(Y_1,0)\right)(Z_1,0)=-\dfrac{H^f(Y_1,Z_1)}{f}(0,X_2)$\\
			\item		$R\left((X_1,0),(0,Y_2)\right)(0,Z_2)=-fg^F(Y_2,Z_2)\left(\nabla^B_{X_1}\grad^B(f),0 \right)$\\
			\item

			$R\left((0,X_2),(0,Y_2)\right)(0,Z_2)=\left( 0,R^F(X_2,Y_2)Z_2\right)+\grad^B f(f)\left(g^F(X_2,Z_2)(0,Y_2)-g^F(Y_2,Z_2)(0,X_2)\right)$
			
		\end{enumerate}
		\quad \\
		where $H^f(U,V)=U(V(f))-(\nabla^B_UV)(f)$ for $U,V\in \chi(F)$ is the Hessian of $f$.

	\end{lem}
	\section{Harmonic vector fields on warped products} In this section, we  determine the harmonicity conditions for vector fields on the warped product $M =B\times_fF$ with  $f: B\longrightarrow ]0,+\infty[$.\\ \text{ }\\
	Let $\{e'_i\}_{i=1_{,\cdots,}m}$ be an orthonomal basis of $(B,g_B)$ and $\{e''_i\}_{i=1_{,\cdots,}n}$  an orthonomal basis of $(F,g_F)$. Then $\{e_i\}_{i=1,_{,\cdots,}m+n}$ is an orthonomal basis of $(M,g)$ with $e_i=(e'_i,0)$ for $i=1,_{\cdots},m$ and\\ $e_{i+m}=\dfrac{1}{f}(0,e''_i)$ for $i=1,_{\cdots},n$. Hence, for $i=1,_{\cdots},m \text{ } \text{ }$  	\begin{eqnarray*}
		\nabla_{e_i}{e_i}&=&(\nabla^B_{e_i'}{e_i'},0)\\ \\
		\nabla_{e_i}V&=&\nabla_{(e_i',0)}{(V_1,0)}+\nabla_{\left(e_i',0\right)}{(0,V_2)}\\&=&\left(\nabla^B_{e_i'}V_1,0\right)+\frac{(e_i',0)(f)}{f}\left(0,V_2\right).
	\end{eqnarray*}
	Moreover
	\begin{eqnarray*}
		(e_i',0)f&=&g\bigg((\grad^Bf,0),(e_i',0)\bigg)=g_B(\grad^Bf,e_i')=e_i'(f).
	\end{eqnarray*}
	Hence, for $i=1,_{\cdots},m \text{ } \text{ }$
	
	\begin{eqnarray*}
		\nabla_{e_i}V=\left(\nabla^B_{e_i'}V_1,0\right)+\dfrac{e_i'(f)}{f}\left(0,V_2\right).
	\end{eqnarray*}
	so
	\begin{eqnarray*}
		\nabla_{e_i}\nabla_{e_i}V&=&\nabla_{(e_i',0)}(\nabla^B_{e_i'}V_1,0)+       \nabla_{(e_i',0)}\left(\dfrac{(e_i',0)(f)}{f}(0,V_2)\right)\\&=&\left(\nabla^B_{e_i'}\nabla^B_{e_i'}V_1,0\right)+\dfrac{e_i'(f)}{f}\nabla_{(e_i',0)}(0,V_2)+(e_i',0)\left(\dfrac{(e_i',0)f}{f}\right)\left(0,V_2\right)\\&=&\left(\nabla^B_{e_i'}\nabla^B_{e_i'}V_1,0\right)+\dfrac{e_i'(f)^2}{f^2}\left(0,V_2\right)+\dfrac{e_i'e_i'(f)}{f}\left(0,V_2\right)-\dfrac{e_i'(f)^2}{f^2}\left(0,V_2\right)\\&=&\left(\nabla^B_{e_i'}\nabla^B_{e_i'}V_1,0\right)+\dfrac{e_i'e_i'(f)}{f}\left(0,V_2\right)
	\end{eqnarray*}
	and
	\begin{eqnarray*}
		\nabla_{\nabla_{e_i}e_i}V&=&\nabla_{\left(\nabla^B_{e_i'}e_i',0\right)}(V_1,0)+\nabla_{\left(\nabla^B_{e_i'}e_i',0\right)}(0,V_2)\\&=&\left(\nabla^B_{\nabla^B_{e_i'}e_i'}V_1,0\right)+\dfrac{\left(\nabla^B_{e_i'}e_i'\right)(f)}{f}(0,V_2).
	\end{eqnarray*}

	For \text{ }$i=m+1,_{\cdots},m+n$, we have
	\begin{eqnarray*}
		\nabla_{e_i}{e_i}&=&\dfrac{1}{f}\nabla_{(0,e_i'')}\dfrac{1}{f}(0,e_i'')\\&=&\dfrac{1}{f}\left(\dfrac{1}{f}\nabla_{(0,e_i'')}(0,e_i'')+(0,e_i'')\big(\dfrac{1}{f}\big)(0,e_i'')\right)\\&=&\dfrac{1}{f^2}\left(\nabla_{(0,e_i'')}(0,e_i'')\right)
		\\&=&\dfrac{1}{f^2}\left(0,\nabla^F_{e_i''}{e_i''}\right)-\dfrac{1}{f}\left(\grad^B f,0\right)	\end{eqnarray*}
	
	because  $$(0,e_i'')(f)=g\left((\grad^B f,0),(0,e_i'')\right)=g_B(\grad^Bf,0)+f^2g_F(0,e_i'')=0.$$\\
	
	For $i=m+1,_{\cdots},m+n$, we compute
	\begin{eqnarray*}
		\nabla_{e_i}V&=&\dfrac{1}{f}\nabla_{(0,e_i'')}(V_1,V_2)\\&=&\dfrac{1}{f}\left( \nabla_{(0,e_i'')}(V_1,0)+\nabla_{(0,e_i'')}(0,V_2)\right)\\&=& \dfrac{1}{f}\left(\left(0,\nabla^F_{e_i''}V_2\right)-fg^F(V_2,e_i'')\left(\grad^Bf,0\right)+\dfrac{(V_1,0)(f)}{f}(0,e_i'')\right)\\&=&\dfrac{1}{f}\left(0,\nabla^F_{e_i''}V_2\right)-g^F(V_2,e_i'')\left(\grad^Bf,0\right)+\dfrac{(V_1,0)(f)}{f^2}(0,e_i'')\\ \nabla_{e_i}V&=& \dfrac{1}{f}\left(0,\nabla^F_{e_i''}V_2\right)-g^F(V_2,e_i'')\left(\grad^Bf,0\right)+\dfrac{V_1(f)}{f^2}(0,e_i'').
	\end{eqnarray*}
	therefore
	\begin{eqnarray*}
		\nabla_{e_i}\nabla_{e_i}V&=&\dfrac{1}{f}\left[   \nabla_{(0,e_i'')}\left(\dfrac{1}{f}(0,\nabla^F_{e_i''}V_2)\right)-\nabla_{(0,e_i'')}\left(g^F(V_2,e_i'')(\grad^Bf,0)\right)\right.\\&&\left.+\nabla_{(0,e_i'')}\left(\dfrac{(V_1,0)f}{f^2}(0,e_i'')\right) \right].
	\end{eqnarray*}
	We now compute the previous terms and sum on $i=m+1,_{\cdots},m+n$:\\ i.)
	\begin{eqnarray}
	\nabla_{(0,e_i'')}\left[\dfrac{1}{f}(0,\nabla^F_{e_i''}V_2)\right]&=&\dfrac{1}{f}\nabla_{(0,e_i'')}\left(0,\nabla^F_{e_i''}V_2\right)   +(0,e_i'')\left(\dfrac{1}{f}\right)\left(0,\nabla^F_{e_i''}V_2\right) \nonumber\\&=&\dfrac{1}{f}\nabla_{(0,e_i'')}\left(0,\nabla^F_{e_i''}V_2\right)\nonumber\\&=&\dfrac{1}{f}\left\{\left(0,\nabla^F_{e_i''}\nabla^F_{e_i''}V_2\right)-fg^F(e_i'',\nabla^F_{e_i''}V_2)\left(\grad^Bf,0\right)\right\}\nonumber\\&=&\dfrac{1}{f}\left(0,\nabla^F_{e_i''}\nabla^F_{e_i''}V_2\right)-g^F(e_i'',\nabla^F_{e_i''}V_2)\left(\grad^Bf,0\right)\nonumber.
	\end{eqnarray}
	ii.)$	\nabla_{(0,e_i'')}\left[g^F(V_2,e_i'')(\grad^Bf,0)\right]\\$
	\begin{eqnarray*}
		&=&g^F(V_2,e_i'')\dfrac{(\grad^Bf,0)(f)}{f}(0,e_i'')
		+(0,e_i'')\left( g^F(V_2,e_i'')\right)\left(\grad^Bf,0\right)\nonumber\\&=&g^F(V_2,e_i'')\dfrac{(\grad^Bf)(f)}{f}(0,e_i'') +e_i''\left(g^F(V_2,e_i'')\right)\left(\grad^Bf,0\right)\nonumber\\&=&\dfrac{\grad^Bf(f)}{f}(0,V_2) +e_i''\left(g^F(V_2,e_i'')\right)\left(\grad^Bf,0\right)\nonumber
	\end{eqnarray*}
	iii.)
	\begin{eqnarray*}
		\nabla_{(0,e_i'')}\left[\dfrac{(V_1,0)(f)}{f^2}(0,e_i'')\right]&=&\dfrac{1}{f^2}\big[\nabla_{(0,e_i'')}(V_1,0)(f)\big](0,e_i'')+(0,e_i'')\left(\dfrac{1}{f^2}(V_1,0)(f)\right)(0,e_i'')\\&=&\dfrac{1}{f^2}(V_1,0)(f)\nabla_{(0,e_i'')}(0,e_i'')\\&=&\dfrac{1}{f^2}V_1(f)\left((0,\nabla^F_{e_i''}e_i'')-f(\grad^Bf,0)\right)\\&=&\dfrac{1}{f^2}V_1(f)(0,\nabla^F_{e_i''}e_i'')-\dfrac{n}{f}V_1(f)(\grad^Bf,0)
	\end{eqnarray*}
	Hence, gathering all the terms, and summing on $i=m+1,_{\cdots},m+n$, we obtain
	
	\begin{eqnarray*}
		\nabla_{e_i}\nabla_{e_i}V&=& \dfrac{1}{f^2}\left(0,\nabla^F_{e_i''}\nabla^F_{e_i''}V_2\right)-\dfrac{1}{f}e_i''(g^F(V_2,e_i''))\left(\grad^Bf,0\right)-\dfrac{\grad^Bf(f)}{f^2}(0,V_2)\\& &+\dfrac{1}{f^3}V_1(f)\left(0,\nabla^F_{e_i''}e_i''\right)-n\dfrac{V_1(f)}{f^2}\left(\grad^Bf,0\right)-\dfrac{1}{f}g^F(e_i'',\nabla^F_{e_i''}V_2)\left(\grad^Bf,0\right), \end{eqnarray*}

	\begin{eqnarray*}
		\nabla_{\nabla_{e_i}e_i}V&=& \nabla_{  \dfrac{1}{f^2}\left(0,\nabla^F_{e_i''}{e_i''}\right)}(V_1,0)-\nabla_{\dfrac{1}{f}(\grad^B f,0)}(V_1,0)+\nabla_{\dfrac{1}{f^2}\left(0,\nabla^F_{e_i''}{e_i''}\right)}(0,V_2)\\&&-\nabla_{\dfrac{1}{f}\left(\grad^Bf,0\right)}(0,V_2)  \\&=& \dfrac{1}{f^3}V_1(f)\left(0,\nabla^F_{e_i''}{e_i''}\right)-\dfrac{n}{f}\left(\nabla^B_{\grad^Bf}V_1,0\right)-n\dfrac{grad^Bf(f)}{f^2}(0,V_2)\\&&+\dfrac{1}{f^2}\left(0,\nabla^F_{\nabla^F_{e_i''}{e_i''}}V_2\right)-\dfrac{1}{f}g^F\left(V_2,\nabla^F_{e_i''}e_i''\right)(\grad^Bf,0).
	\end{eqnarray*}

	Hence, summing on the index $i$,
	\begin{eqnarray*}
		\nabla^*\nabla V&=&\nabla_{\nabla_{e_i}e_i}V-\nabla_{e_i}\nabla_{e_i}V\\&=& \dfrac{1}{f^3}V_1(f)\left(0,\nabla^F_{e_i''}{e_i''}\right)-\dfrac{n}{f}\left(\nabla^B_{\grad^Bf}V_1,0\right)- n\dfrac{ \grad^Bf(f)}{f^2}(0,V_2)\\&&+\dfrac{1}{f^2}\left(0,\nabla^F_{\nabla^F_{e_i''}{e_i''}}V_2\right)-\dfrac{1}{f}g^F(V_2,\nabla^F_{e_i''}e_i'')(\grad^Bf,0)+ \left(\nabla^B_{\nabla^B_{e_i'}e_i'}V_1,0\right)\\&&+\dfrac{\left(\nabla^B_{e_i'}e_i'\right)(f)}{f}(0,V_2)
		-(\nabla^B_{e_i'}\nabla^B_{e_i'}V_1,0)-\dfrac{e_i'e_i'(f)}{f}(0,V_2)-\dfrac{1}{f^2}\left(0,\nabla^F_{e_i''}\nabla^F_{e_i''}V_2\right)\\&&+\dfrac{1}{f}e_i''(g^F(V_2,e_i''))\left(\grad^Bf,0\right)+\dfrac{ \grad^Bf(f)}{f^2}(0,V_2)-\dfrac{1}{f^3}V_1(f)\left(0,\nabla^F_{e_i''}e_i''\right)\\&&+n\dfrac{V_1(f)}{f^2}\left(\grad^Bf,0\right)+\dfrac{1}{f}g^F(e_i'',\nabla^F_{e_i''}V_2)\left(\grad^Bf,0\right) 
		\\&=&
		\left( \nabla{^*}\nabla V_1-\dfrac{n}{f}\nabla^B_{\grad^Bf}V_1+\dfrac{2}{f}g^F(e_i'',\nabla^F_{e_i''}V_2)\grad^Bf+n\dfrac{V_1(f)}{f^2}\grad^Bf\textbf{,\quad}\dfrac{1}{f^2}\nabla{^*}\nabla V_2\right.\\&&\left.-\dfrac{e_i'e_i'(f)}{f}V_2+\dfrac{\left(\nabla^B_{e_i'}{e_i'}\right)(f)}{f}V_2+(1-n)\dfrac{\grad f(f)}{f^2}(0,V_2)\right) .
	\end{eqnarray*}
	\begin{lem}
		Let $(B,g_B)$ and $(F, g_F)$ be Riemannian manifolds and $f: B\longrightarrow \R^*_+$ a smooth function on B. Let $\{e'_i\}_{i=1_{,\cdots,}m}$ be an orthonomal basis of $(B,g_B)$ and $\{e''_i\}_{i=1_{,\cdots,}n}$  an orthonomal basis of $(F,g_F)$. Then   a vector field $V=V_1+V_2$  on $M=B\times_fF$ is a harmonic vector field if and only if \begin{align}\label{Eq: for2}\begin{cases}\displaystyle 
		\nabla{^*}\nabla V_1-\dfrac{n}{f}\nabla^B_{\grad^Bf}V_1+\dfrac{2}{f}(\di^F V_2)\grad^Bf+n\dfrac{V_1(f)}{f^2}\grad^Bf=0 ,\\ \\\displaystyle\dfrac{1}{f^2}\nabla{^*}\nabla V_2+\dfrac{\Delta^B(f)}{f}V_2+(1-n)\dfrac{\grad f(f)}{f^2}V_2=0 ,
		\end{cases}\end{align}
		\\where $\Delta^B(f)=-\tr H^{f}$.	\end{lem}
	
	\begin{prop}
		Let $G$ be a 2-dimensional Riemannian Lie group equipped with a left-invariant metric, $f:I\longrightarrow \R^*_+$ a smooth function on $\R$, $V_1=\phi(t)\partial_t$ a vector field on $\R$ and $V_2$ a unit vector field on $G$. Then $V=\phi(t)\partial_t+V_2$ is a harmonic vector field on the warped product $I\times_fG$ if  $f(t)= \sqrt{2\kappa_0 t^2+c_1t+c_2} $ on $I$ such that 
		\begin{itemize}
			\item$I=]-\infty,-\frac{c_2}{c_1}]$ if $\kappa_0=0$ and $c_1>0$
			\item $I=[-\frac{c_2}{c_1},+\infty[$ if $\kappa_0=0$ and $c_1<0$
			\item $I=\R$ if $c_1^2-8c_2\kappa_0<0$ and $\kappa_0>0$
			\item $I=[t_1,t_2[$ if $c_1^2-8c_2\kappa_0\geq0$ and $\kappa_0>0$
			\item $I=]-\infty,t_1[\cup]t_2,+\infty[$ if $c_1^2-8c_2\kappa_0\geq0$ and $\kappa_0<0$
		\end{itemize} and $\phi$ is a solution on $I$  of the differential equation \begin{eqnarray}\label{Eq: fr10} (E_0): 
		x''+3\frac{f'}{f}x'-3\left(\frac{f'}{f}\right)^2x+2\kappa_1\frac{f'}{f}=0\end{eqnarray}
		where $\displaystyle \kappa_0=\big<\nabla^*\nabla V_2,V_2>, \kappa_1=\sum_{i=1}^{n}g^F(e'_i,\nabla_{e'_i}V_2)$, $t_1=\frac{-c_1-\sqrt{c_1^2-8c_2\kappa_0}}{4\kappa_0}$ and $t_2=\frac{-c_1+\sqrt{c_1^2-8c_2\kappa_0}}{4\kappa_0}$ with $t_1\leq t_2$.
	\end{prop}
	Note that the $ODE$ $(E_0)$ always admits a solutions defined on the whole of $I$.
	\begin{lem}\label{Eq: fr}\cite{Amann}
		Consider the non-homogenous linear differential equation\begin{eqnarray}\label{Eq: rew1}
		\sum_{k=0}^{n}h_k(t)y^k(t)=\si(t)
		\end{eqnarray}If the functions $h_k (k\in \{0,1,2,\cdots,n\})$ are continuous on the same interval $ I $, with $h_n$ not zero on this interval, there exist $ n $ linearly independent solutions of the homogeneous differential equation associated with (\ref {Eq: rew1}) on $I$, and the most general solution is the sum of a linear combination with arbitrary constant coefficients of these n functions and a particular solution of (\ref {Eq: rew1}).
	\end{lem}

	Take $(B,g_B)$ to be $(\mathbb{R},dt^2)$
	and $V=V_1+V_2$ on $\R\times_fF$ where $V_1=\phi(t)\partial_t$ is a vector field on $\R$ and  $\displaystyle V_2=\sum_{i=1}^{n}a_ie_i$  a vector field on $F$, where the $a_i's$ are functions on $F$. Note that for  an orthonomal basis $\{e_i\}_{i=1_{,\cdots,}n}$  of $(F,g)$, the Levi-Civita connection $\nabla$ of $(F,g)$ and $\{\partial_t\}$ the canonical vector field on $\mathbb{R}$, we have:
	\begin{eqnarray*}
		\nabla_{e_i}V_2&=&\sum_{j=1}^{n}\bigg(a_j\nabla_{e_i}e_j+e_i(a_j)e_j\bigg) \qquad \text{for}\quad i=1_{,\cdots,}n,\\
		\nabla_{e_i}\nabla_{e_i}V_2&=&\sum_{j=1}^{n}\bigg(a_j\nabla_{e_i}\nabla_{e_i}e_j+e_ie_i(a_j)e_j+2e_i(a_j)\nabla_{e_i}e_j\bigg) \qquad \text{for}\quad i=1_{,\cdots,}n,\end{eqnarray*}

	\begin{eqnarray*}
		\nabla_{\nabla_{e_i}e_i}V_2&=&\sum_{j=1}^{n}\bigg(a_j\nabla_{\nabla_{e_i}e_i}e_j+(\nabla_{e_i}e_i)(a_j)e_j\bigg) \qquad \text{for}\quad i=1_{,\cdots,}n.
	\end{eqnarray*}

	We now compute the previous terms and sum on $i=1,\cdots,n$ to obtain\\
	\begin{align*}
	\nabla^*\nabla V_2=& \sum_{i,j=1}^{n}\bigg(a_j\nabla_{\nabla_{e_i}e_i}e_j+(\Delta a_j)e_j-2e_i(a_j)\nabla_{e_i}e_j-a_j\nabla_{e_i}\nabla_{e_i}e_j\bigg).
	\end{align*}
	We also have, $$g(e_i,\nabla_{e_i}V_2)=e_i(a_i)+\sum_{j=1}^{n}\bigg(a_jg(e_i,\nabla_{e_i}e_j)\bigg),$$ 
	\\
	
	and\begin{eqnarray*}
		\nabla{^*}\nabla V_1&=&\nabla_{\nabla_{\partial_t}\partial_t}\varphi(t)\partial_t-\nabla_{\partial_t }\nabla_{\partial_t}\phi(t)\partial_t=-\phi''(t)\partial_t\\
		\grad^Bf&=&f'(t)\partial_t\\
		\nabla^B_{\grad^Bf}V_1&=&f'(t)\phi'(t)\partial_t\\ e_i'(e_i'(f))&=&\partial_t(\partial_t(f))=f''(t).\end{eqnarray*} 
	Then V is a harmonic vector field if and only if \\

	\begin{align*}\label{Eq: for2}\begin{cases}\displaystyle
	\phi''(t)+n\frac{f'(t)}{f(t)}\phi'(t)-n\left(\frac{f'(t)}{f(t)}\right)^2\phi-2\left(\sum_{i=1}^{n}\left(e_i(a_i)+\sum_{j=1}^{n}a_jg(e_i,\nabla_{e_i}e_j)\right)\right)\frac{f'(t)}{f(t)}&=0,
	\\\displaystyle \sum_{j=1}^{n}\left((\Delta a_j)e_j+\sum_{i=1}^{n}\left(\left(a_j\nabla_{\nabla_{e_i}e_i}e_j-2e_i(a_j)\nabla_{e_i}e_j-a_j\nabla_{e_i}\nabla_{e_i}e_j\right)\right)\right)\\-\displaystyle\bigg(f(t)f''(t)+2f'(t)^2\bigg)\sum_{i=
		1}^{n}a_je_j=0.\end{cases}
	\end{align*}
	
	\section{Harmonic vector fields on $\R\times_fG$: $G$ unimodular Lie group}
	
	\subsection{Vector fields constructed from unit left-invariant vector fields}
	
	In this section, we assume that $F$ is a three-dimensional connected Riemannian Lie group $G$ equipped with an invariant metric. We determine   harmonic vector fields on $(\mathbb{R}\times_fG,g)$ with $f:\mathbb{R}\longrightarrow ]0;+\infty[$. Let $\{e_1,e_2,e_3\}$ be an orthonormal basis on $G$.
	
	Let $V_2=ae_1+be_2+ce_3$ such that $a^2+b^2+c^2=1$, $V_1=\phi(t)\partial_t$ and consider the vector field $V=V_1+V_2$ on $\R\times_fG$.\\ 
	
	\begin{prop}[\cite{Milnor76}]
		
		Let $G$ be a three-dimensional unimodular connected Lie group, $\mathfrak{g}$ its Lie algebra and $g$ a left-invariant Riemannian metric on $G$. Then there exists an orthonormal basis $\{e_1,e_2,e_3\}$ of $\mathfrak{g}$ such that \begin{eqnarray}\label{Eq: Phare4}
		[e_2, e_3]=\lambda_1e_1,\quad [e_3, e_1]=\lambda_2e_2,\quad [e_1, e_2]=\lambda_3e_3\end{eqnarray} where $\lambda_1, \lambda_2,\lambda_3$ are constants. \\
		\begin{center}\textbf{Table 1}: Three-dimensional unimodular Lie groups\\ \text{ }\\
			\begin{tabular}{|c|c|}
				\hline &\\
				Signs of $\lambda_1,\lambda_2,\lambda_3$\quad	&\quad \text{ }\quad \quad \quad \quad Associated Lie Groups   \\  &\\\hline
				+,+,+& \quad\quad \quad$SU(2)$ or $SO(3)$\\
				\hline +,+,-&\quad\quad \quad $SL(2,\R)$ or $O(1,2)$\\\hline +,+,0&\quad\quad \quad $\mathbb{E}(2)$\\\hline+,0,-&\quad\quad \quad $\mathbb{E}(1,1)$ \\\hline +,0,0&\quad\quad \quad $\mathbb{H}^3$\\\hline0,0,0&\quad\quad \quad$\mathbb{R}\oplus\mathbb{R}\oplus\mathbb{R}$\\
				\hline
			\end{tabular}
		\end{center}\text{ }
	\end{prop}

	Then the Levi-Civita connection $\nabla$ is determined by \cite{Milnor76}:

	\begin{equation*}
	\nabla_{e_1}e_1=0, \quad\nabla_{e_1}e_2=\mu_1e_3,\quad \nabla_{e_1}e_3=-\mu_1e_2;
	\end{equation*} 
	\begin{equation*}
	\nabla_{e_2}e_1=-\mu_2e_3,\quad \nabla_{e_2}e_2=0,\quad\nabla_{e_2}e_3=\mu_2e_1;
	\end{equation*}
	\begin{equation*}
	\nabla_{e_3}e_1=\mu_3e_2,\quad \nabla_{e_3}e_2=-\mu_3e_1, \quad\nabla_{e_3}e_3=0;
	\end{equation*}
	where $$\mu_i=\frac{1}{2}(\lambda_1+\lambda_2+\lambda_3)-\lambda_i, \text{ }i=1,2,3.$$ 
	
	From this, we have
	\begin{equation*}
	\nabla_{\nabla_{e_i}e_i}e_j=\nabla_{e_i}\nabla_{e_i}e_i=g(e_i,\nabla_{e_i}e_j)=0
	\end{equation*}
	\begin{equation*}
	\nabla_{e_i}\nabla_{e_i}e_j=-\mu_i^2e_j, \quad i,j=1,2,3.
	\end{equation*}
	Hence, $V$ \text{ }\text{is a harmonic vector field if and only if }\text{ }
	\begin{eqnarray}\label{Eq: fr1}
	\begin{cases}
	\phi''(t)+3\frac{f'(t}{f(t)}\phi'(t)-3\left(\frac{f'(t)}{f(t)}\right)^2\phi(t)=0,
	\\a\left(\mu_2^2+\mu_3^2-f(t)f''(t)-2f'(t)^2\right)=0,\\ b\left(\mu_1^2+\mu_3^2-f(t)f''(t)-2f'(t)^2\right)=0,\\c\left(\mu_2^2+\mu_1^2-f(t)f''(t)-2f'(t)^2\right)=0.
	\end{cases}.\end{eqnarray}
	
	The first equation of (\ref{Eq: fr1}) always admits a non-trivial solution.\\ For the other equations \begin{eqnarray}\label{Eq: fr4}
	\begin{cases}
	a\left(\mu_2^2+\mu_3^2-f(t)f''(t)-2f'(t)^2\right)=0,\\ b\big(\mu_1^2+\mu_3^2-f(t)f''(t)-2f'(t)^2\big)=0,\\c\big(\mu_2^2+\mu_1^2-f(t)f''(t)-2f'(t)^2\big)=0,
	\end{cases}.\end{eqnarray} we distinguish six cases according to the classification of three-dimensional unimodular Lie groups given in Table 1 and the domain $I$ of the solutions of the non-linear  differential equation $yy''+2y'^2=\epsilon, \epsilon\in \R$.\\ \\ First, recall that the non-linear differential equation $yy''+2y'^2=\epsilon, \epsilon\in \R$ always admits solutions by \rm{Cauchy-Peano-Arzelà}:

	\begin{thm}[Cauchy-Peano-Arzelà]\cite{Amann}
		Let $ f:U=I\times V\longrightarrow E$ be a continuous function, where $I $ and $V$ are open sets of $\R$ and $E$ (a Banach space) respectively, and let $(t_0, x_0)\in U $. Then there exists a neighbourhood $I$ of $t_0$ and an application $X: J\longrightarrow V$ such that $\forall t \in J, X'(t)=f(t, X (t)); X (t_0)=x_0 $.
		More precisely, if $T=]t_0-\al, t_0+\al [\times B(x_0,r)$ has its adherence in $U$ and if $ M\al <r$ with
		$M=\sup{\|f(t, x)\|, (t, x) \in T}$, then there exists
		a solution to the Cauchy's problem $X'(t)=f(t,X (t))$ with $ X(t_0)=x_0$, defined on $J=]t_0-\al, t_0+\al[$.
		This neighborhood $T$ of $(t_0,x_0)$ is sometimes called the "safety barrel".
	\end{thm}
	\begin{enumerate}
		\item $G=\R\oplus\R\oplus\R$. In this case, $\lambda_1=\lambda_2=\lambda_3=0$ so that $\mu_1=\mu_2=\mu_3=0$, and (\ref{Eq: fr4}) becomes \begin{eqnarray*}
			\begin{cases}a\left(f(t)f''(t)+2f'(t)^2\right)=0,\\b\left(f(t)f''(t)+2f'(t)^2\right)=0,\\c\left(f(t)f''(t)+2f'(t)^2\right)=0.
			\end{cases}.\end{eqnarray*} Since $(a,b,c)\ne(0,0,0)$, we have $f(t)f''(t)+2f'(t)^2=0$, so that $f(t)=\left(3e^{c_1}t+c_2\right)^{\frac{1}{3}}$, $c_1,c_2\in \R$ and $\phi$ is defined on $I=]-\frac{c_2}{3e^{c_1}},+\infty[$. We obtain a harmonic vector field $V=\phi(t)\partial_t+V_2$ on $I\times_f G$ for $f(t)=\left(3e^{c_1}t+c_2\right)^{\frac{1}{3}}$ and $V_2=ae_1+be_2+ce_3$ for any $a,b,c\in \R, a^2+b^2+c^2=1$.
		\\\item $G=\mathbb{H}^3$ is the Heisenberg group. In this case, $\lambda_1>0,\lambda_2=0, \lambda_3=0$, we have $\mu_{1}^2=\mu_{2}^2=\mu_{3}^2=\mu^2>0$, and (\ref{Eq: fr4}) becomes \begin{eqnarray*}
			\begin{cases}a\big(2\mu^2-f(t)f''(t)-2f'(t)^2\big)=0,\\b\big(2\mu^2-f(t)f''(t)-2f'(t)^2\big)=0,\\c\big(2\mu^2-f(t)f''(t)-2f'(t)^2\big)=0.
		\end{cases}\end{eqnarray*}Since $(a,b,c)\ne(0,0,0)$ and $\mu\ne0$, we necessarily have $f(t)f''(t)+2f'(t)^2=2\mu^2$ and obtain a harmonic vector field $V=\phi(t)\partial_t+V_2$ on $I\times_f G$ for $f$ solution  of  $y y''+2y'^2=2\mu^2$ on $I$ and $V_2=ae_1+be_2+ce_3$, for any $a,b,c\in \R, a^2+b^2+c^2=1$.
		
		\item $G=\mathbb{E}(1,1)$. In this case, $\lambda_1>0,\lambda_2=0,\lambda_3<0$, we have\\ $\mu_{3}=-\mu_1=\frac{1}{2}(\lambda_1-\lambda_3)>0, \mu_2=\frac{1}{2}(\lambda_1+\lambda_3)$ hence $\mu_{1}^2=\mu_3^2\ne0$, and the System (\ref{Eq: fr4}) becomes
		\begin{eqnarray*}
			\begin{cases}
				a\big(\mu_1^2+\mu_2^2-f(t)f''(t)-2f'(t)^2\big)=0,\\ 
				b\big(2\mu_1^2-f(t)f''(t)-2f'(t)^2\big)=0,\\
				c\big(\mu_1^2+\mu_2^2-f(t)f''(t)-2f'(t)^2\big)=0.
			\end{cases}	
		\end{eqnarray*} 
		If 
		\begin{enumerate}
			\item $b=0$ and  $(a,c)\ne (0,0)$, we have $f(t)f''(t)+2f'(t)^2=\mu_1^2+\mu_2^2$.
			We obtain a harmonic vector field $V=\phi(t)\partial_t+V_2$ on $I\times_f G$ for  $f$ solution of $y y''+2y'^2=\mu_1^2+\mu_2^2$ on $I$ and $V_2=ae_1+ce_3$ for any $a,c\in \R, a^2+c^2=1$.
			\item $b=\pm1$ and  $a=c=0$, we have $f(t)f''(t)+2f'(t)^2=2\mu_1^2$.
			We obtain a harmonic vector field $V=\phi(t)\partial_t+V_2$ on $I\times_f G$ for  $f$ solution of $y y''+2y'^2=2\mu_1^2$ on $I$ and $V_2=\pm e_2$.
		\end{enumerate}
		
		\item $G=\mathbb{E}(2)$. In this case, $\lambda_1>0,\lambda_2>0, \lambda_3=0$ and we have\\ 
		$\mu_{2}=-\mu_1=\frac{1}{2}(\lambda_1-\lambda_2), \mu_3=\frac{1}{2}(\lambda_1+\lambda_2)>0$ hence $\mu_{1}^2=\mu_2^2$, and (\ref{Eq: fr4})  becomes
		\begin{eqnarray*}
			\left\lbrace \begin{array}{lll}
				a\big(\mu_1^2+\mu_3^2-f(t)f''(t)-2f'(t)^2\big)=0,\\ 
				b\big(\mu_1^2+\mu_3^2-f(t)f''(t)-2f'(t)^2\big)=0,\\
				c\big(2\mu_1^2-f(t)f''(t)-2f'(t)^2\big)=0.
			\end{array}\right.	
		\end{eqnarray*} 
		
		If
		\begin{enumerate}
			\item $c=0$ and $ (a,b)\ne (0,0)$, we have $f(t)f''(t)+2f'(t)^2=\mu_1^2+\mu_3^2$. 
			We obtain a harmonic vector field $V=\phi(t)\partial_t+V_2$ on $I\times_f G$ for  $f$ solution of $y y''+2y'^2=\mu_1^2+\mu_3^2$ on $I$, $V_2=ae_1+be_2$ for any $a,b\in \R,\, a^2+b^2=1$.
			
			\item $c=\pm1$ and  $a=b=0$, we have $f(t)f''(t)+2f'(t)^2=2\mu_1^2$.
			We obtain a harmonic vector field $V=\phi(t)\partial_t+V_2$ on $I\times_f G$ for  $f$ solution of $y y''+2y'^2=2\mu_1^2$ on $I$ and $V_2=\pm e_3$.
		\end{enumerate}
		
		\item $G=SL(2,\R)$ or $O(1,2)$. In this case $\lambda_1>0,\lambda_2>0,\lambda_3<0$. We distinguish two cases:
		\begin{enumerate}
			\item $\lambda_1=\lambda_2>0>\lambda_3$, we have $\mu_1=\mu_2=\frac{\lambda_3}{2}<0$, $\mu_3=\frac{1}{2}(2\lambda_1-\lambda_3)$ and (\ref{Eq: fr4})  becomes
			\begin{eqnarray*}
				\begin{cases}
					a\big(\mu_1^2+\mu_3^2-f(t)f''(t)-2f'(t)^2\big)=0,\\ 
					b\big(\mu_1^2+\mu_3^2-f(t)f''(t)-2f'(t)^2\big)=0,\\
					c\big(2\mu_1^2-f(t)f''(t)-2f'(t)^2\big)=0.
				\end{cases}	
			\end{eqnarray*} 
			The different cases are:
			\begin{enumerate}
				\item $c=0$ and $ (a,b)\ne (0,0)$. We have $f(t)f''(t)+2f'(t)^2=\mu_1^2+\mu_3^2$. 
				We obtain a harmonic vector field $V=\phi(t)\partial_t+V_2$ on $I\times_f G$ for  $f$ solution of $y y''+2y'^2=\mu_1^2+\mu_3^2$ on $I$, $V_2=ae_1+be_2$ for any $a,b\in \R, a^2+b^2=1$.
				
				\item $c=\pm1$ and  $a=b=0$. We have $f(t)f''(t)+2f'(t)^2=2\mu_1^2$.
				We obtain a harmonic vector field $V=\phi(t)\partial_t+V_2$ on $I\times_f G$ for  $f$ solution of $y y''+2y'^2=2\mu_1^2$ on $I$ and $V_2=\pm e_3$.
			\end{enumerate}
			
			\item $\lambda_1>\lambda_2>0>\lambda_3$. From (\ref{Eq: fr4}), we have
			\begin{enumerate}
				\item  $a\ne0$, $b=c=0$ and $f(t)f''(t)+2f'(t)^2=\mu_3^2+\mu_2^2$ .We obtain a harmonic vector field $V=\phi(t)\partial_t+V_2$ on $I\times_f G$ for $f$ solution of $y y''+2y'^2=\mu_3^2+\mu_2^2$ on $I$ and $V_2=\pm e_1$.			
				\item  $b\ne0$, $a=c=0$, and $f(t)f''(t)+2f'(t)^2=\mu_3^2+\mu_1^2$.
				We obtain a harmonic vector field $V=\phi(t)\partial_t+V_2$ on $I\times_f G$ for $f$ solution of $y y''+2y'^2=\mu_1^2+\mu_3^2$ on $I$ and $V_2=\pm e_2$.
				\item $c\ne0$, $a=0=b$ and  $f(t)f''(t)+2f'(t)^2=\mu_3^2+\mu_2^2$.
				We obtain a harmonic vector field $V=\phi(t)\partial_t+V_2$ on $I\times_f G$ for $f$ solution of $y y''+2y'^2=\mu_2^2+\mu_2^2$ on $I$ and  $V_2=\pm e_3$.				
			\end{enumerate}
		\end{enumerate}
		\item $G=SU(2)$ or $SO(3)$. In this case $\lambda_1>0,\lambda_2>0,\lambda_3>0$. We distinguish three subcases:
		\begin{enumerate}
			\item $\lambda_1=\lambda_2=\lambda_3>0$ so that $\mu_1=\mu_2=\mu_3=\mu=\frac{1}{2}\lambda_1$. From (\ref{Eq: fr4}), since $(a,b,c)\ne (0,0,0)$ we have $f(t)f''(t)+2f'(t)^2=2\mu^2$.
			We obtain a harmonic vector field $V=\phi(t)\partial_t+V_2$ on $I\times_f G$ for $f$ solution of $y y''+2y'^2=2\mu^2$ on $I$ and $V_2=ae_1+be_2+ce_3$ for any $b,c\in \R, a^2+b^2+c^2=1$.

			\item $\lambda_1>\lambda_2=\lambda_3$: we have
			\begin{enumerate}
				\item $a=0$, $ (b,c)\ne (0,0)$ and $f(t)f''(t)+2f'(t)^2=\mu_1^2+\mu_2^2$. 
				We obtain a harmonic vector field $V=\phi(t)\partial_t+V_2$ on $I\times_f G$ for $f$ solution of $y y''+2y'^2=\mu_1^2+\mu_2^2$ on $I$ and  $V_2=be_2+ce_3$ for any $b,c\in \R, b^2+c^2=1$.
				
				\item $a=\pm1$, $b=c=0$ and  $f(t)f''(t)+2f'(t)^2=2\mu_2^2$.
				We obtain a harmonic vector field $V=\phi(t)\partial_t+V_2$ on $I\times_f G$ for $f$ solution of $y y''+2y'^2=2\mu_2^2$ on $I$ and $V_2=\pm e_1$.
			\end{enumerate}
			
			\item $\lambda_1=\lambda_2>\lambda_3$: We have  
			\begin{enumerate}
				\item $c=0$, $ (a,b)\ne (0,0)$ and  $f(t)f''(t)+2f'(t)^2=\mu_1^2+\mu_3^2$. 
				We obtain a harmonic vector field $V=\phi(t)\partial_t+V_2$ on $I\times_f G$ for $f$ solution of $y y''+2y'^2=2\mu_1^2+\mu_3^2$ on $I$ and $V_2=ae_1+be_2$ for any $b,a\in \R, b^2+a^2=1$.
				
				\item $c=\pm1$, $a=b=0$ and $f(t)f''(t)+2f'(t)^2=2\mu_1^2$.
				We obtain a harmonic vector field $V=\phi(t)\partial_t+V_2$ on $I\times_f G$ for $f$ solution of $y y''+2y'^2=2\mu_1^2$ on $I$ and $V_2=\pm e_3$.
			\end{enumerate}
			\item The case $\lambda_1>\lambda_2>\lambda_3$ is similar to the second case of $SL(2,\R)$.\\
		\end{enumerate}
	\end{enumerate}
	\begin{prop}
		Let $G$ be a connected  three-dimensional unimodular Riemannian Lie group  equipped with an invariant metric and $V=V_1+V_2$  a vector field on the warped product $I\times_fG$, $f: I\longrightarrow ]0,+\infty[$ with $V_1=\phi(t)\partial_t$ and $V_2=ae_1+be_2+ce_3$ a unit left-invariant vector field on $G$. Let $\{e_1,e_2,e_3\}$ be an orthonormal basis of the Lie algebra satisfying (\ref{Eq: Phare4}) and $\lambda_1,\lambda_2,\lambda_3$ the structure constants. Then $V$ is a harmonic vector field if and only if $f$  is a solution  on $I\subset \R$ of the ODE $$(E_1): xx''+2x'^2=\epsilon,\quad\epsilon\in \R$$, 
		$\phi$ is a solution, on $I$, of the ODE \begin{eqnarray}\label{Eq: fr7} (E_0): 
		\displaystyle  x''+3\dfrac{f'}{f}x'-3\left(\dfrac{f'}{f}\right)^2x=0\end{eqnarray} and one of the following cases occurs:
		\begin{enumerate}
			\item $\lambda_1=\lambda_2=\lambda_3=0$, $\displaystyle f(t)=\left(3e^{c_1}t+c_2\right)^{\frac{1}{3}}$, $I=]-\frac{c_2}{3e^{c_1}},+\infty[$, and $V_2=ae_1+be_2+ce_3$ for any $a,b,c\in \R, a^2+b^2+c^2=1$,
			\item $\lambda_1>0,\lambda_2=\lambda_3=0,\epsilon=2\mu_1^2$ and  $V_2=ae_1+be_2+ce_3$ for any $a,b,c\in \R, a^2+b^2+c^2=1$,
			\item  $\lambda_1>0,\lambda_2=0,\lambda_3<0,$
			\begin{enumerate}
				\item $\epsilon=\mu_1^2+\mu_2^2$ and  $V_2=ae_1+ce_3$ for any $a,c\in \R, a^2+c^2=1$,
				\item $\epsilon=2\mu_1^2$ and $V_2=\pm e_2$,
			\end{enumerate}
			\item  $\lambda_1>0,\lambda_2>0,\lambda_3=0$,
			\begin{enumerate}
				\item $\epsilon=\mu_1^2+\mu_3^2$ and  $V_2=ae_1+be_2$ for any $a,b\in \R, a^2+b^2=1$,
				\item $\epsilon=2\mu_1^2$ and  $V_2=\pm e_3$,
			\end{enumerate}
			
			\item  $\lambda_1=\lambda_2>0>\lambda_3$
			\begin{enumerate}
				\item $\epsilon=\mu_1^2+\mu_3^2$ and  $V_2=ae_1+be_2$ for any $a,b\in \R, a^2+b^2=1$,
				\item $\epsilon=2\mu_1^2$ and $V_2=\pm e_3$,
			\end{enumerate}
			\item  $\lambda_1>\lambda_2>0>\lambda_3$ or $\lambda_1>\lambda_2>\lambda_3,\epsilon=\mu_i^2+\mu_j^2$ and  $V_2=\pm e_k$ $i\ne j\ne k$,
			
			\item 
			$\lambda_1=\lambda_2=\lambda_3>0,\epsilon=2\mu_1^2$ and  $V_2=ae_1+be_2+ce_3$ for any $a,b,c\in \R, a^2+b^2+c^2=1$,
			\item $\lambda_1>\lambda_2=\lambda_3$,\begin{enumerate}
				\item $\epsilon=\mu_1^2+\mu_2^2$ and  $V_2=be_2+ce_3$ for any $b,c\in \R, c^2+b^2=1$,
				\item $\epsilon=2\mu_2^2$ and $V_2=\pm e_1$,
			\end{enumerate}
			
			\item $\lambda_1=\lambda_2>\lambda_3>0$,
			\begin{enumerate}
				\item $\epsilon=\mu_1^2+\mu_3^2$ and   $V_2=ae_1+be_2$ for any $a,b\in \R, a^2+b^2=1$
				\item $\epsilon=2\mu_1^2$ and  $V_2=\pm e_3$.
			\end{enumerate}
			
		\end{enumerate}
		\end{prop}
	\begin{rem}
		Note that if $f$ is a non-zero positive constant on $\R$  and $G$ is a three-dimensional unimodular connected  Lie group  equipped with an invariant metric, then $V=(V_1,V_2)$ is harmonic vector field on $\R\times_fG$ if and only if $\phi$ is an affine function, $\mu_j=\mu_i\ne \mu_k$ and $V_2=\pm e_k$, $i,j,k\in\{1,2,3\}$   or $\mu_1=\mu_2=\mu_3=0$ and  $V_2$ any vector field on $G$. This recovers the results in \cite{Calvaruso}. 
	\end{rem}
	
	\subsection{Vector fields constructed from non left-invariant vector fields}
	
	In this subsection, we construct new examples of harmonic vector fields from non left-invariant vector fields on a unimodular Lie group.
	
	\begin{exam}Let $G=\R\oplus\R\oplus\R$ and consider $V_2=a(x,y,z)\frac{\partial}{\partial x}$. Using Relation (\ref{Eq: for2}),  $V=\varphi(t)\partial_t+V_2$ is a harmonic vector field  if and only if   \begin{eqnarray*}\begin{cases}
				\phi''(t)+3\frac{f'(t)}{f(t)}\phi'(t)-3\left(\frac{f'(t)}{f(t)}\right)^2\phi=2k\frac{f'(t)}{f(t)},\\ \\a_x=k,\quad k\in \R , \\\\f(t)f''(t)+2f'(t)^2=\epsilon \\ \\\De a=\epsilon a.
			\end{cases}	
		\end{eqnarray*} 
		This implies that:\\ \textbf{Case 1}: $\epsilon>0$, then \begin{eqnarray*}
			\begin{cases}\displaystyle a(x,y,z)=\cos(v_1y+v_2z)+\sin(v_1y+v_2z),\quad v_1^2+v_2^2=\epsilon, \quad \\ \\\displaystyle f(t)f''(t)+2f'(t)^2=\epsilon, \quad t\in I\subset \R ,\\\\ \displaystyle
				\phi''(t)+3\frac{f'(t)}{f(t)}\phi'(t)-3\left(\frac{f'(t)}{f(t)}\right)^2\phi=2k\frac{f'(t)}{f(t)}. \end{cases}	
		\end{eqnarray*} 
		\textbf{Case 2}: $\epsilon<0$, then 
		
		\begin{eqnarray*}
			\begin{cases}\displaystyle a(x,y,z)=\exp(v_1y+v_2z),\quad v_1^2+v_2^2=-\epsilon,\quad \\ \\\displaystyle f(t)f''(t)+2f'(t)^2=\epsilon, \quad t\in I\subset \R ,\\\\\displaystyle
				\phi''(t)+3\frac{f'(t)}{f(t)}\phi'(t)-3\left(\frac{f'(t)}{f(t)}\right)^2\phi=2k\frac{f'(t)}{f(t)}. \end{cases}	
		\end{eqnarray*} 
		\\ \textbf{Case 3:} $\epsilon=0$, then 	
		\begin{eqnarray*}
			\begin{cases} \displaystyle a(x,y,z)=kx+b(y,z) \quad\text{$b$ is a harmonic function on $\R^2$},\\ \\\displaystyle f(t)=\left(3e^{c_1}x+c_2\right)^{\frac{1}{3}} \quad\text{and}\quad  I=]-\frac{c_2}{3e^{c_1}},+\infty[,
				
				\quad\\ \\\displaystyle 
				\phi''(t)+3\frac{f'(t)}{f(t)}\phi'(t)-3\left(\frac{f'(t)}{f(t)}\right)^2\phi=2k\frac{f'(t)}{f(t)}.\end{cases}		
		\end{eqnarray*}
	\end{exam}

	\begin{prop}\label{Eq: P1}
		The vector field $V=\varphi(t)\partial_t+V_2$ is a harmonic vector field on $I\times_f\big(\R\oplus\R\oplus\R\big)$ if and only if 
		\begin{eqnarray*}
			\begin{cases}
				f(t)f''(t)+2f'(t)^2=\epsilon, \quad t\in I\subset \R ,\quad \epsilon\ne0 ,\\ \\\phi''(t)+3\frac{f'(t)}{f(t)}\phi'(t)-3\left(\frac{f'(t)}{f(t)}\right)^2\phi=2\kappa\frac{f'(t)}{f(t)},		\end{cases}	
		\end{eqnarray*} with $a(x,y,z)=\cos(v_1y+v_2z)+\sin(v_1y+v_2z), v_1^2+v_2^2=-\epsilon$ if $\epsilon<0$ and  $a(x,y,z)=\exp(v_1y+v_2z)
		,v_1^2+v_2^2=\epsilon$ if $\epsilon>0$, or 
		
		\begin{eqnarray*}
			\begin{cases}\displaystyle a(x,y,z)=\kappa_1x+b(y,z) ,\\\displaystyle f(t)=\left(3e^{c_1}t+c_2\right)^{\frac{1}{3}} \quad\text{and}\quad  t\in I=]-\frac{c_2}{3e^{c_1}},+\infty[ ,
				
				\quad\\ \displaystyle
				\phi''(t)+3\frac{f'(t)}{f(t)}\phi'(t)-3\left(\frac{f'(t)}{f(t)}\right)^2\phi=2k\frac{f'(t)}{f(t)},\end{cases}		
		\end{eqnarray*}

		where $b(y,z)$ is a harmonic function on $\R^2$ and $\kappa, \kappa_1\in \R$.	
		
	\end{prop}
	\begin{exam}Let $G=\mathbb{H}^3$ be the Heisenberg group of real $3\times 3$ upper-triangular matrices of the form
		\begin{equation}
		\begin{array}{c}
		A
		\end{array}
		=\left (
		\begin{array}{cccc}
		1 &x &y\\
		0 &1 &z\\
		0 & 0 & 1
		\end{array}\right ),
		\end{equation}
		endowed with the left-invariant metric given by $dx^2+\big(dy-xdz\big)^{2}+(dz)^{2}$. We identify $\mathbb{H}^3$ with ${\mathbb{R}}^{3}$, endowed with this metric. The left-invariant vector fields
		\begin{equation}
		\begin{array}{lcr}
		e_{1}=\frac{\partial}{\partial x}, &e_{2}=\frac{\partial}{\partial y}, &e_{3}=\frac{\partial}{\partial z}+x\frac{\partial}{\partial y},
		\end{array}
		\end{equation}
		constitute an orthonormal basis of the Lie algebra $\mathfrak{g}$ of $ \mathbb{H}^3$ and the corresponding Levi-Civita connection is determined by
		\begin{equation}\label{Eq: for5}
		\begin{array}{lcr}
		\nabla_{e_{1}}e_{2}=\nabla_{e_{2}}e_{1}=-\frac{1}{2}e_{3}, \\ \nabla_{e_{1}}e_{3}=-\nabla_{e_{3}}e_{1}=\frac{1}{2}e_{2}, \\ \nabla_{e_{2}}e_{3}=\nabla_{e_{3}}e_{2}=\frac{1}{2}e_{1},
		\end{array}
		\end{equation}
		where the remaining covariant derivatives vanish.
		\\ By(\ref{Eq: for5}), we have $\mu_1=-\frac{1}{2}=-\mu_2=-\mu_3$. Take $V_2=a(x,y,z)e_1$, then the vector field $V=\varphi(t)\partial_t+V_2$ is a  harmonic vector field if \begin{eqnarray*}
			\begin{cases}
				\phi''(t)+3\frac{f'(t)}{f(t)}\phi'(t)-3\left(\frac{f'(t)}{f(t)}\right)^2\phi=2\kappa\frac{f'(t)}{f(t)}, \\\\a_x=\kappa,\quad \kappa\in \R,\\\\ a_y=0,\\\\  a_z+xa_y=0,\\ \\f(t)f''(t)+2f'(t)^2=\epsilon, \quad\epsilon\in \R,\\\\\De a=(\epsilon-\frac{1}{2})a.
			\end{cases}
		\end{eqnarray*} 
		Therefore  
		\begin{eqnarray*}
			\begin{cases} a(x,y,z)=\kappa x+\kappa'  \quad \kappa,\kappa'\in \R,\\ \\f(t)f''(t)+2f'(t)^2=\frac{1}{2}, \quad t\in I\subset \R,\\ \\ \phi''(t)+3\frac{f'(t)}{f(t)}\phi'(t)-3\left(\frac{f'(t)}{f(t)}\right)^2\phi=2\kappa\frac{f'(t)}{f(t)}. \end{cases}	
		\end{eqnarray*} 
	\end{exam}

	\begin{prop}\label{Eq: P2}
		Let $V_1=\phi(t)\partial_t$  on $\R$ and $V_2=a(x,y,z)e_1$  on $\mathbb{H}^3$, then $V=V_1+V_2$ is a harmonic vector field on $I\times_f \mathbb{H}^3$ if and only if 
		\begin{eqnarray*}
			\begin{cases} a(x,y,z)=\kappa x+\kappa'  \quad \kappa,\kappa'\in \R,\\ \\f(t)f''(t)+2f'(t)^2=\frac{1}{2}, \quad t\in I\subset \R,\\ \\ \phi''(t)+3\frac{f'(t)}{f(t)}\phi'(t)-3\left(\frac{f'(t)}{f(t)}\right)^2=2\kappa\frac{f'(t)}{f(t)}. \end{cases}	
		\end{eqnarray*}
	\end{prop}

	\section{Harmonic vector fields on $\R\times_fG$: $G$ non-unimodular Lie groups}
	\subsection{Vector fields constructed from unit left-invariant vector fields}

	In this subsection, $F$ is now a connected three-dimensional non-unimodular Riemannian Lie group $G$ equipped with an invariant metric and we determine   harmonic vector fields on $(\mathbb{R}\times_fG,g)$ with $f:\mathbb{R}\longrightarrow ]0;+\infty[$. Let $\{e_1,e_2,e_3\}$ be an orthonormal basis on $G$, $V_2=ae_1+be_2+ce_3$ such that $a^2+b^2+c^2=1$ and $V_1=\phi(t)\partial_t$ and consider the vector field $V=V_1+V_2$ on $\R\times_fG$.
	\begin{prop}[\cite{Milnor76}] Let $G$ be a connected three-dimensional Riemannian non-unimodular Lie group, $\mathfrak{g}$ its Lie algebra and $g$ a left-invariant Riemannian metric on $G$. Then there exists an orthonormal basis $\{e_1,e_2,e_3\}$ of $\mathfrak{g}$ such that \begin{eqnarray}\label{Eq: Phare5}
		[e_1, e_2]=\alpha e_2+\beta e_3,\quad [e_1, e_3]=-\beta e_2+\delta e_3, \quad[e_2, e_3]=0.\end{eqnarray} where $\alpha+\delta>0$ and $\al\geq\de$ are constants.
	\end{prop} Then the Levi-Civita connection $\nabla$ is determined by \cite{Milnor76}

	\begin{center}$\nabla_{e_1}e_1=0, \nabla_{e_1}e_2=\beta e_3,\quad \nabla_{e_1}e_3=-\beta e_2$;
	\end{center}
	\begin{center}$\nabla_{e_2}e_1=-\alpha e_2,\quad \nabla_{e_2}e_2=\alpha e_1,\quad\nabla_{e_2}e_3=0$;\end{center} 
	\begin{center}$\nabla_{e_3}e_1=-\delta e_3,\quad \nabla_{e_3}e_2=0,\quad \nabla_{e_3}e_3=\delta e_1$;\end{center} 
	From this, we have 
	\begin{equation*}
	\nabla_{\nabla_{e_1}e_1}e_j=\nabla_{e_2}\nabla_{e_2}e_1=\nabla_{e_3}\nabla_{e_3}e_1=0,\quad\nabla_{e_1}\nabla_{e_1}e_1=\nabla_{e_3}\nabla_{e_3}e_2=\nabla_{e_2}\nabla_{e_2}e_3=0  
	\end{equation*}	
	\begin{equation*}
	\nabla_{\nabla_{e_3}e_3}e_2=\de\be e_3,\quad\nabla_{e_3}\nabla_{e_3}e_3=-\de\be e_2,\quad
	\nabla_{\nabla_{e_2}e_2}e_3=-\al\be e_2,\quad\nabla_{e_2}\nabla_{e_2}e_2=\al\be e_3   ,
	\end{equation*}
	\begin{equation*}
	\nabla_{e_1}\nabla_{e_1}e_2=-\be^2e_2,\quad \nabla_{e_1}\nabla_{e_1}e_3=-\be^2e_3,\quad \nabla_{e_2}\nabla_{e_2}e_1=-\al^2e_1,\quad \nabla_{e_2}\nabla_{e_2}e_2=-\al^2e_2,
	\end{equation*}
	\begin{equation*}
	\nabla_{e_3}\nabla_{e_3}e_1=-\de^2e_1,\quad \nabla_{e_3}\nabla_{e_3}e_3=-\de^2e_3,\quad g(e_i,\nabla_{e_i}e_i)=g(e_1,\nabla_{e_1}e_i)=0 ,
	\end{equation*}
	\begin{equation*} 
	g(e_3,\nabla_{e_3}e_2)=g(e_2,\nabla_{e_2}e_3)=0,\quad
	g(e_2,\nabla_{e_2}e_1)=-\al,\quad g(e_3,\nabla_{e_3}e_1)=-\de.
	\end{equation*}	
	Hence $V$ \text{ }is a harmonic vector field if and only if 
	
	\begin{eqnarray}\label{Eq: fr2} \begin{cases}
	\phi''(t)+3\frac{f'(t)}{f(t)}\phi'(t)-3\frac{f'(t)^2}{f(t)^2}\phi(t)=-2a(\de+\al)\frac{f'(t)}{f(t)},\\	\\a\left(\alpha^2+\delta^2-f(t)f''(t)-2f'(t)^2\right)=0,\\ \\
	b\left(\alpha^2+\beta^2-f(t)f''(t)-2f'(t)^2\right)-\beta(\alpha+\delta)c=0,\\ \\
	c\left(\delta^2+\beta^2-f(t)f''(t)-2f'(t)^2\right)+\beta(\alpha+\delta)b=0.		
	\end{cases}\end{eqnarray}
	The first equation of (\ref{Eq: fr2}) has a solution by Theorem (\ref{Eq: fr}) and the remaining system becomes: 
	
	\begin{eqnarray*} \begin{cases}a\left(\alpha^2+\delta^2-f(t)f''(t)-2f'(t)^2\right)=0,\\ \\
			b\left(\alpha^2+\beta^2-f(t)f''(t)-2f'(t)^2\right)-\beta(\alpha+\delta)c=0,\\ \\
			c\left(\delta^2+\beta^2-f(t)f''(t)-2f'(t)^2\right)+\beta(\alpha+\delta)b=0.		
	\end{cases}\end{eqnarray*}We analyse the different cases:
	\begin{enumerate}
		\item $\al=\de>0$: \begin{enumerate}
			\item Suppose $a\ne0$, then $f(t)f''(t)+2f'(t)^2=2\al^2>0$  and $c=b=0$.
			We obtain a harmonic vector field $V=(\phi(t)\partial_t,\pm e_1)$ on $I\times_f G$ for  $f$ solution of $y y''+2y'^2=2\al^2$ on $I$.
			\item Suppose $a=0$: Necessarily
			\begin{enumerate}
				\item  $b=0$ and $c=\pm1$ hence $\be=0$ and $f(t)f''(t)+2f'(t)^2=2\al^2$.
				We obtain a harmonic vector field $V=(\phi(t)\partial_t,\pm e_3)$ on $I\times_f G$ for  $f$ solution of $y y''+2y'^2=2\al^2$ on $I$.
				
				\item $c=0$ and $b=\pm1$ hence $\be=0$ and $f(t)f''(t)+2f'(t)^2=2\al^2$.
				We obtain a harmonic vector field $V=(\phi(t)\partial_t,\pm e_2)$ on $I\times_f G$ for  $f$ solution of $y y''+2y'^2=2\al^2$ on $I$.
				\item $b\ne0,c\ne0, \be=0$, and $f(t)f''(t)+2f'(t)^2=\al^2+\be^2$. 
				We obtain a harmonic vector field $V=(\phi(t)\partial_t,V_2)$ on $I\times_f G$ for  $f$ solution of $y y''+2y'^2=\al^2+\be^2$ on $I$  and $V_2=be_2+ce_3$ for any $b,c\in \R;  b^2+c^2=1$.
				
			\end{enumerate}
		\end{enumerate}
		\item $\al>\de>0$:
		\begin{enumerate}
			\item Suppose $a\ne0$, then $f(t)f''(t)+2f'(t)^2=\al^2+\de^2$ and $\begin{cases}
			b(\be^2-\de^2)=\be(\al+\de)c,\\c(\be^2-\al^2)=-\be(\al+\de)b.	
			\end{cases}$\\
			We obtain a harmonic vector field $V=(\phi(t)\partial_t,V_2)$ on $I\times_f G$ for  $f$ solution of $y y''+2y'^2=\al^2+\de^2$ on $I$ and $V_2=ae_1+be_2+ce_3$ for any $a,b,c\in \R, a^2+b^2+c^2=1, a\ne0$.
			\item Suppose that $a=0$:
			\begin{enumerate}
				\item $b=0, c=\pm1$ hence $\beta=0$, we have $f(t)f''(t)+2f'(t)^2=\de^2$.  
				We obtain a harmonic vector field $V=(\phi(t)\partial_t,\pm e_3)$ on $I\times_f G$ for  $f$ solution of $y y''+2y'^2=\de^2$ on $I$.
				\item $ c=0, b=\pm1$ hence $\beta=0$, we have $f(t)f''(t)+2f'(t)^2=\al^2$.  
				We obtain a harmonic vector field $V=(\phi(t)\partial_t,\pm e_2)$ on $I\times_f G$ for  $f$ solution of $y y''+2y'^2=\al^2$ on $I$.
				\item $b\ne0,c\ne0$, we have $\be=bc(\al-\de)$ and $f(t)f''(t)+2f'(t)^2=\be^2+b^2\al^2+c^2\de^2$. We obtain a harmonic vector field $V=(\phi(t)\partial_t,V_2)$ on $I\times_f G$ for $f$ solution of $y y''+2y'^2=\be^2+b^2\al^2+c^2\de^2$ on $I$ and $V_2=be_2+ce_3$ for any $a,b,c\in \R, b^2+c^2=1, b\ne 0,c\ne0$.
			\end{enumerate}
		\end{enumerate}
		\item $\al>\de=0$: 
		\begin{enumerate}
			\item $a=0$:
			\begin{enumerate}
				\item $b=0, c=\pm1$ and $\beta=0$, we have $f(t)f''(t)+2f'(t)^2=0$, \\so that  $\displaystyle {f(t)=\left(3e^{c_1}t+c_2\right)^{\frac{1}{3}}}$ and $\phi$ is defined on $I=]-\frac{c_2}{3e^{c_1}},+\infty[$.
				We obtain a harmonic vector field $V=(\phi(t)\partial_t,\pm e_3)$ on $I\times_f G$ for  $f(t)=\left(3e^{c_1}t+c_2\right)^{\frac{1}{3}}$.
				
				\item $ c=0, b=\pm1$ hence $\beta=0$, we have $f(t)f''(t)+2f'(t)^2=\al^2$ 
				We obtain a harmonic vector field $V=(\phi(t)\partial_t,\pm e_2)$ on $I\times_f G$ for  $f$ solution of $y y''+2y'^2=\al^2$ on $I$.
				\item $b\ne0,c\ne0$ we have $\be=bc\al$ and $f(t)f''(t)+2f'(t)^2=\be^2+b^2\al^2$. We obtain a harmonic vector field $V=(\phi(t)\partial_t,V_2)$ on $I\times_f G$ for  $f$ solution of $y y''+2y'^2=\be^2+b^2\al^2$ on $I$ and $V_2=be_2+ce_3, b\ne0,c\ne0$.
			\end{enumerate}
			\item $a\ne0$, we have $f(t)f''(t)+2f'(t)^2=\al^2$ and 
			$\begin{cases}\be(b\be-c\al)=0,\\c(\be^2-\al^2)+\be\al b=0.		
			\end{cases}$\\  
		
			This implies $c=0$ and $\be=0$ or  $b=c=0$.
			We obtain a harmonic vector field $V=(\phi(t)\partial_t,V_2)$ on $I\times_f G$ for $f$ solution of $y y''+2y'^2=\al^2$ on $I$ and $V_2=ae_1+be_2, a\ne0$, $a^2+b^2=1$.
		\end{enumerate}
		\item $\al>0>\de$:
		\begin{enumerate}
			\item $a=0$, we have
			\begin{eqnarray*}\begin{cases}b(\al^2+\be^2-f(t)f''(t))-\be c(\al+\de)=0,\\		c(\de^2+\be^2-f(t)f''(t))+\be b(\al+\de)=0.
			\end{cases}\end{eqnarray*}
			Then
			\begin{enumerate}
				\item  $b=0$, $c=\pm1$, $\be=0$ and $f(t)f''(t)+2f'(t)^2=\de^2$.			
				We obtain a harmonic vector field $V=(\phi(t)\partial_t,\pm e_3)$ on $I\times_f G$ for for  $f$ solution of $y y''+2y'^2=\de^2$ on $I$.				
				\item  $c=0$, $b=\pm1$, $\be=0$ and $f(t)f''(t)+2f'(t)^2=\al^2$.	We obtain a harmonic vector field $V=(\phi(t)\partial_t,\pm e_3)$ on $I\times_f G$ for for  $f$ solution of $y y''+2y'^2=\al^2$ on $I$.
				
				\item $ b\ne 0,c\ne0$, we have $\be=bc(\al-\de)$ then $f(t)f''(t)+2f'(t)^2=b^2\al^2+c^2\de^2+\be^2$.
				We obtain a harmonic vector field $V=(\phi(t)\partial_t,V_2)$ on $I\times_f G$ for $f$ solution of $y y''+2y'^2=b^2\al^2+c^2\de^2+\be^2$ on $I$ and $V_2=be_2+ce_3$, $(b^2+c^2=1)$.

			\end{enumerate}
			\item $a\ne 0$, we have $f(t)f''(t)+2f'(t)^2=\al^2+\de^2$ and $\begin{cases}
			b(\be^2-\de^2)=\be(\al+\de)c,\\c(\be^2-\al^2)=-\be(\al+\de)b .\end{cases}$
			
			We obtain a harmonic vector field $V=(\phi(t)\partial_t,V_2)$ on $I\times_f G$  for $f$ solution of\\ $y y''+2y'^2=\al^2+\de^2$ on $I$ and  $V_2=ae_1+be_2+ce_3$, with $a^2+b^2+c^2=1, a\ne0.$
		\end{enumerate}
	\end{enumerate}

	\begin{prop}Let $G$ be a connected  three-dimensional Riemannian  non-unimodular  Lie group  equipped with an  invariant metric and $V=V_1+V_2$   a vector field on the warped product $I\times_fG$, $f: I\longrightarrow ]0,+\infty[$ with $V_1=\phi(t)\partial_t$ and $V_2$ a unit left-invariant vector field on $G$.  Let $\{e_1,e_2,e_3\}$ be an orthonormal basis of the Lie algebra satisfying Equation(\ref{Eq: Phare5}) and $\al,\be,\de$ the structure constants. Then $V$ is a harmonic vector field if and only if $f$  satisfies, on $I\subset \R$, the ordinary  differential equation $$(E_1): x x''+2x'^2=\epsilon,\quad\epsilon\in \R$$,  $\phi$ is a solution, on $I$, of the equation \begin{eqnarray}\label{Eq: fr0} (E_a): 
		x''+3\frac{f'}{f}x'-3\left(\frac{f'}{f}\right)^2x=-2a(\de+\al)\frac{f'}{f}\end{eqnarray} and $V_2$ is determined by one of the following conditions:

		\begin{enumerate}
			\item $\al=\de>0$
			\begin{enumerate}
				\item $\be=0$, $\epsilon=2\al^2$ and $V_2=\pm e_1$,
				\item $\be=0$, $\epsilon=2\al^2$ and $V_2=\pm e_2$ or $V_2=\pm e_3$,
				\item $\be=0$,  $\epsilon=\al^2+\be^2$ and $V_2=be_2+ce_3, b^2+c^2=1, b\ne 0,c\ne0$,
			\end{enumerate}
			\item $\al>\de>0$ or $\al>0>\de$
			\begin{enumerate}
				\item $\epsilon=\al^2+\de^2$,  $V_2=ae_1+be_2+ce_3$,$a\ne0$,  $a^2+b^2+c^2=1$ and $\begin{cases}
				b(\be^2-\de^2)=\be(\al+\de)c,\\c(\be^2-\al^2)=-\be(\al+\de)b	 .
				\end{cases}$
				\item $\be=0$, $\epsilon=\de^2, V_2=\pm e_3$,
				\item $\be=0$, $\epsilon=\al^2, V_2=\pm e_2$,
				\item $\be=bc(\al-\de), \epsilon=\be^2+b^2\al^2+c^2\de^2$ and $V_2=be_2+ce_3, c\ne 0,b\ne 0, b^2+c^2=1$.
			\end{enumerate}
			\item $\al>\de=0$
			\begin{enumerate}
				\item $\be=0$, in this case $f(t)=\left(3e^{c_1}t+c_2\right)^{\frac{1}{3}}$, $I=]-\frac{c_2}{3e^{c_1}},+\infty[$ and $V_2=\pm e_3$,
				\item $\be=0$, $\epsilon=\al^2$ and $V_2=\pm e_2$,
				\item $\be=bc\al$, $\epsilon=\be^2+b^2\al^2$ and $V_2=be_2+ce_3$, $b\ne0,c\ne0, b^2+c^2=1$,
				\item $\epsilon=\al^2$ and $V_2=\pm e_1$,
				\item $\be=0$, $\epsilon=\al^2$ and $V_2=ae_1+be_2$, $a\ne 0, a^2+b^2=1$.
			\end{enumerate}

		\end{enumerate}
		\begin{rem}
			Note that if $f$ is a non-zero positive constant on $\R$ and $G$ is a three-dimensional non-unimodular  Lie group  equipped with an invariant metric, then $V=(V_1,V_2)$ is harmonic vector field if and only if $\phi$ is an affine  function on $\R, \be=\de=0$ $V_2=be_2+ce_3$ and $\begin{cases}
			b(\be^2-\de^2)=\be(\al+\de)c,\\c(\be^2-\al^2)=-\be(\al+\de)b.	
			\end{cases}$
		\end{rem}
			
	\end{prop}
	
	\subsection{Vector fields constructed from non left-invariant vector fields}

	In this subsection, we give examples of harmonic vector fields on $\R\times_fG$  constructed from non left-invariant vector fields on the non-unimodular Lie group $G$

	\begin{exam}Let $F$ be $(\R\times H^{2}, g)$ where 
		$H^{2}=\Large\{(x,y)\in {\mathbb{R}}^{2}: y>0\Large\}$  denotes the Poincaré half-plane with Gaussian curvature equal to $-\al$ $(\al>0)$ and $g$  the left-invariant metric given by
		$$g=\frac{1}{\al y^{2}}\Big(dx^2+dy^{2}\Big)+dz^2.$$The left-invariant vector fields
		\begin{equation*}
		\begin{array}{lcr}
		e_{1}=y\sqrt{\al}\frac{\partial}{\partial y}, &e_{2}=y\sqrt{\al}\frac{\partial}{\partial x}, &e_{3}=\frac{\partial}{\partial z},
		\end{array}
		\end{equation*} 
		constitute an orthonormal basis of the Lie algebra $\mathfrak{g}$ of $\R\times H^{2},$ and $$[e_1,e_2]=\sqrt{\al} e_2, \quad[e_3,e_1]=0,\quad[e_3,e_2]=0.$$The corresponding Levi-Civita connection is determined by \begin{equation*}\label{Eq: f5}
		\begin{array}{lcr}
		\nabla_{e_{2}}e_{1}=-\sqrt{\al}e_2,\quad\nabla_{e_{2}}e_{2}=\sqrt{\al} e_1,
		\end{array}
		\end{equation*}
		where the remaining covariant derivatives  vanish.  Take $V_2=b(y)+c(y)$ and use (\ref{Eq: for2}) to see that $V=\varphi(t)\partial_t+V_2$ is a harmonic vector field if and only if 
		\begin{eqnarray*}
			\begin{cases}
				\phi''(t)+3\frac{f'(t)}{f(t)}\phi'(t)-3\left(\frac{f'(t)}{f(t)}\right)^2\phi=0,\\\\\Delta b=\left(f(t)f"(t)+2f'(t)^2-\al\right)b,\\\\\Delta c=(f(t)f"(t)+2f'(t)^2)c.
		\end{cases}\end{eqnarray*}

		This is equivalenty to
		
		\begin{eqnarray*}
			\begin{cases}
				\phi''(t)+3\frac{f'(t)}{f(t)}\phi'(t)-3\left(\frac{f'(t)}{f(t)}\right)^2\phi=0,\\ \\\Delta b=\left(f(t)f"(t)+2f'(t)^2-\al\right)b,\\\\\Delta c=(f(t)f"(t)+2f'(t)^2)c,\\\\f(t)f"(t)+2f'(t)^2=\epsilon.
		\end{cases}\end{eqnarray*}Hence 
		
		\begin{eqnarray*}
			\begin{cases}
				\phi''(t)+3\frac{f'(t)}{f(t)}\phi'(t)-3\left(\frac{f'(t)}{f(t)}\right)^2\phi=0,\\\\ y^2b''(y)=(\epsilon-\al)b(y), \\\\y^2c''(y)=\al c(y),\\\\f(t)f"(t)+2f'(t)^2=\epsilon.
		\end{cases}\end{eqnarray*} We have 
		
		\begin{eqnarray}\label{Eq: r1}
		\begin{cases}\displaystyle c(y)=y^{\frac{1}{2}}\bigg(\kappa_1+\kappa_2\ln(y)\bigg) \text{ }\text{if}\text{ } \al=4\epsilon,\\

		\displaystyle	c(y)=\kappa_1y^{\frac{1}{2}+\frac{1}{2}\sqrt{1-4\al}}+\kappa_2 y^{\frac{1}{2}-\frac{1}{2}\sqrt{1-4\al}} \text{ }\text{ }\text{ }\text{if}\text{ } 4\epsilon<\al ,
		
		\\ \displaystyle c(y)=y^{\frac{1}{2}}\bigg(\kappa_1\cos(y\sqrt{4\al-1})+\kappa_2\sin(y\sqrt{4\al-1})\bigg) \text{ }\text{ }\text{ }\text{if}\text{ } \al<4\epsilon,
		\end{cases}\end{eqnarray} and

		\begin{eqnarray}\label{Eq: r2}
		\begin{cases}\displaystyle b(y)=y^{\frac{1}{2}}\bigg(\kappa_1+\kappa_2\ln(y)\bigg) \text{ }\text{if}\text{ } 4\epsilon=5\al ,\\

		\displaystyle b(y)=\kappa_1y^{\frac{1}{2}+\frac{1}{2}\sqrt{1-4\al}}+\kappa_2y^{\frac{1}{2}-\frac{1}{2}\sqrt{1-4\al}} \text{ }\text{ }\text{ }\text{if}\text{ } 5\al<4\epsilon ,
		
		\\ \displaystyle b(y)=y^{\frac{1}{2}}\bigg(\kappa_1\cos(y\sqrt{4\al-1})+\kappa_2\sin(y\sqrt{4\al-1})\bigg) \text{ }\text{ }\text{ }\text{if}\text{ } 4\epsilon<5\al,
		\end{cases}\end{eqnarray} with $ \kappa_1,\kappa_2\in \R$.
	\end{exam}

	\begin{prop}\label{Eq: P3}
		Let $V_1=\phi(t)\partial_t$ be a vector field on $\R$ and $V_2=a(ye_1+b(y)e_2+c(y)e_3$ be vector fields on $(\R\times H^{2}, g)$ where 
		$H^{2}=\Large\{(x,y)\in {\mathbb{R}}^{2}: y>0\Large\}$  denotes the Poincaré half-plane with Gaussian curvature equal to $-\al(\al>0)$. Then $V=V_1+V_2$ is a harmonic vector field on $I\times_fG$ if and only if 
		\begin{eqnarray*}
			\begin{cases}
				\phi''(t)+3\frac{f'(t)}{f(t)}\phi'(t)-3\left(\frac{f'(t)}{f(t)}\right)^2\phi=0,\quad t\in I,\\f(t)f"(t)+2f'(t)^2=\epsilon, \text{ }\epsilon \in \R,
			\end{cases}\end{eqnarray*} and $V_2=b(y)e_2+c(y)e_3$ with $b, c$ defined in Equations (\ref{Eq: r1}) and (\ref{Eq: r2}).
		
	\end{prop}
	
	\section{Harmonic maps on warped product}
	In this section, we determine the horizontal part of the tension field on the  warped product $B\times_fF$ (cf 2.7) and  study the existence of  vector fields on $\R\times_f G$ which are harmonic maps, where $G$ is a three-dimensional  Lie group  equipped with an invariant metric. To calculate $S(V)$ where $V=V_1+V_2$, we write 
	\begin{eqnarray*}
		S(V)=S_1(V)+S_2(V) \quad \text{where }\quad S_1(V)=\sum_{i=1}^{m}R(\nabla_{e_i}V,V)e_i \quad \text{and}\quad S_2(V)=\sum_{i=m+1}^{m+n}R(\nabla_{e_i}V,V)e_i.
	\end{eqnarray*}Then
	
	\begin{eqnarray*}
		S_1(V)&=& \sum_{i=1}^{m}R\left( \nabla_{(e_i',0)}V,(V_1,V_2)\right)(e_i',0)\\&=&\sum_{i=1}^{m}R\left( (\nabla^B_{e_i'}V_1,0),(V_1,0)\right)(e_i',0)+R\left( (\nabla^B_{e_i'}V_1,0),(0,V_2)\right)(e_i',0)+\\&&\frac{e_i'(f)}{f}R\left((0,V_2),(V_1,0)\right)(e_i',0)+\frac{e_i'(f)}{f}R\left((0,V_2),(0,V_2)\right)(e_i',0)\\&=&	\sum_{i=1}^{m}
		\left( R^B(\nabla^B_{e_i'}V_1,e_i')V_1,0\right)-\frac{e_i'(f)}{f^2}H^f(V_1,e_i')(0,V_2)+\dfrac{1}{f}H^f(\nabla^B_{e_i'}V_1,e_i')(0,V_2)
	\end{eqnarray*}

	\begin{eqnarray*}\displaystyle
		S_2(V)&=&\displaystyle\sum_{i=1}^{n}\frac{1}{f}R\left( \nabla_{(0,e_i'')}V,V\right)(0,e_i'')\\\displaystyle&=&\sum_{i=1}^{n}	\frac{1}{f^2}R\left((0,\nabla^F_{e_i''}V_2),(V_1,0)\right)(0,e_i'')+\frac{1}{f^2}R\left((0,\nabla^F_{e_i''}V_2),(0,V_2)\right)(0,e_i'')-\\&&\displaystyle\frac{1}{f}g^F(V_2,e_i'') R\left((\grad f,0),(V_1,0)\right)(0,e_i'')-\frac{1}{f}g^F(V_2,e_i'')  R\left((\grad f,0),(0,V_2)\right)(0,e_i'')
		\\&&+\frac{V_1(f)}{f^3}R\left((0,e_i''),(0,V_2)\right)(0,e_i'')+\frac{V_1(f)}{f^3}R\left((0,e_i''),(V_1,0)\right)(0,e_i'')\\\displaystyle&=&\sum_{i=1}^{n}\left(\frac{1}{f^2}\left(0,R^F(\nabla_{e_i''}V_2,V_2)e_i'' \right)+\frac{V_1(f)}{f^3}(0,R^F(e_i'',V_2)e_i'')\right)+\|V_2\|^2(\nabla^B_{\grad f}\grad f,0)+\\\displaystyle&&\frac{\grad f(f)}{f^2}\bigg(\di(V_2)(0,V_2)-\big(0,\nabla^F_{V_2}V_2\big)\bigg)+n\dfrac{V_1(f)}{f^2}\big(\nabla^B_{V_1}\grad f,0\big)+\frac{1}{f}\di(V_2)(\nabla^B_{V_1}\grad  f,0)\\\displaystyle&&+\frac{V_1(f)}{f^3}\grad f(f)(n-1)(0,V_2)
		\\
		\displaystyle&=&\frac{1}{f^2}(0,S(V_2))+\frac{\grad f(f)}{f^2}\bigg(\di(V_2)(0,V_2)-\big(0,\nabla^F_{V_2}V_2\big)\bigg)+n\dfrac{V_1(f)}{f^2}\big(\nabla^B_{V_1}\grad f,0\big)\\\displaystyle&&+\frac{V_1(f)}{f^3}\sum_{i=1}^{n}(0,R^F(e_i'',V_2)e_i'')+\|V_2\|^2(\nabla^B_{\grad f}\grad f,0)+\frac{1}{f}\di(V_2)(\nabla^B_{V_1}\grad  f,0)\\\displaystyle&&+\frac{V_1(f)}{f^3}\grad f(f)(n-1)(0,V_2)
	\end{eqnarray*}
	Hence 
	
	\begin{eqnarray*}\displaystyle
		S(V)&=&S_1(V)+S_2(V)\\\displaystyle&=&\left( (S(V_1),0)-\sum_{i=1}^{m}\bigg(\frac{e_i'(f)}{f^2}H^f(V_1,e_i')-\dfrac{1}{f}H^f(\nabla^B_{e_i'}V_1,e_i')\bigg)(0,V_2)+n\dfrac{V_1(f)}{f^2}\big(\nabla^B_{V_1}\grad f,0\big)+\right.\\\displaystyle&&\left.\frac{1}{f^2}(0,S(V_2))+\frac{\grad f(f)}{f^2}\bigg(\di(V_2)(0,V_2)-\big(0,\nabla^F_{V_2}V_2\big)\bigg)+\frac{V_1(f)}{f^3}\grad f(f)(n-1)(0,V_2)\right.\\\displaystyle&&\left.\frac{V_1(f)}{f^3}\sum_{i=1}^{n}(0,R^F(e_i'',V_2)e_i'')+\|V_2\|^2(\nabla^B_{\grad f}\grad f,0)+\frac{1}{f}\di(V_2)(\nabla^B_{V_1}\grad  f,0)\right)\\
		\displaystyle&=&\left(S(V_1)+n\dfrac{V_1(f)}{f^2}\nabla^B_{V_1}\grad f+\|V_2\|^2\nabla^B_{\grad f}\grad f+\frac{1}{f}\di(V_2)\nabla^B_{V_1}\grad f; \frac{1}{f^2}S(V_2)\right.\\\displaystyle&&\left.-\sum_{i=1}^{m}\bigg(\frac{e_i'(f)}{f^2}H^f(V_1,e_i')-\dfrac{1}{f}H^f(\nabla^B_{e_i'}V_1,e_i')\bigg)V_2+\frac{\grad f(f)}{f^2}\bigg(\di(V_2)V_2-\nabla^F_{V_2}V_2\bigg)+\right.\\\displaystyle&&\left.\frac{V_1(f)}{f^3}\sum_{i=1}^{n}R^F(e_i'',V_2)e_i''+\frac{V_1(f)}{f^3}\grad f(f)(n-1)V_2\right).
	\end{eqnarray*} Take  $B=\R$, $V_1=\varphi(t)\partial_t$ on $\R$ and $V_2$  a vector field on $(F,g)$, then

	\begin{eqnarray*}
		S(V)&=&\left( n\dfrac{\phi^2f'f''}{f^2}\partial_t+\|V_2\|^2f'f''\partial_t+\dfrac{1}{f}\di(V_2)\phi f''\partial_t; \frac{1}{f^2}S(V_2)-\frac{f'\phi f''}{f^2}V_2+\frac{\phi' f''}{f}V_2+	\right.\\\displaystyle&&\left.\frac{f'^2}{f^2}\di(V_2)V_2
		-\frac{f'^2}{f^2}\nabla_{V_2}{V_2}+ \frac{\phi f'^3}{f^3}(n-1)V_2+\frac{\phi f'}{f^3}\sum_{i=1}^{n}R(e_i'',V_2)e_i''\right).
	\end{eqnarray*}
	
	Remark that if $f$ is a non-zero positive constant  on $\R$, then a harmonic vector field $V=V_1+V_2$ on the warped product $\R \times _f F$ is  a harmonic map if and only if  $V_2$ is a harmonic map on $F$.

	\begin{prop}[Unimodular groups]
		Let $G$ be a connected three-dimensional Riemannian unimodular  Lie group  equipped with an invariant metric and $V=V_1+V_2$ a harmonic vector field on the warped product $I\times_fG, I\subset \R$, $f: \R\longrightarrow ]0,+\infty[$ with $V_1=\phi(t)\partial_t$ and $V_2=ae_1+be_2+ce_3$ a unit left-invariant vector field on $G$. Then 
		$$\displaystyle \nabla_{V_2}{V_2}=bc(\mu_2-\mu_3)e_1+ac(\mu_3-\mu_1)e_2+ab(\mu_1-\mu_2)e_3$$
		 and 
		$$\displaystyle S(V_2)=A_1bce_1+A_2ace_2+A_3abe_3$$ where $\displaystyle A_1=\mu_2^2(\mu_3-\mu_1)+\mu_3^2(\mu_1-\mu_2)$,  
		$\displaystyle A_2=\mu_1^2(\mu_2-\mu_3)+\mu_3^2(\mu_1-\mu_2)$,
		$\displaystyle A_3=\mu_1^2(\mu_2-\mu_3)+\mu_2^2(\mu_3-\mu_1)$ and  $V$ is harmonic map on $I\times_f G$ if and only if one of the following cases occurs:
		\begin{enumerate}
			\item $G=\R\oplus\R\oplus\R: \lambda_1=\lambda_2=\lambda_3=0, f(t)=\be>0, I=\R, V_2=ae_1+be_2+ce_3$ for any $a,b,c\in \R$, $a^2+b^2+c^2=1$ and $\phi(t)=\ga_1 t+\ga_2, \ga_1,\ga_2\in \R$.
			
			\item $\displaystyle G=\mathbb{H}^3: \lambda_1>0, \lambda_2=\lambda_3=0, f(t)=\varepsilon\mu_1t+\be \text{ }(\varepsilon=\pm 1), \be\in\R$ and $\displaystyle  I=]-\frac{\be}{\mu_1},+\infty[$ for $\varepsilon=1$, $\displaystyle  I=]-\infty,\frac{\be}{\mu_1}[$ for $\varepsilon=-1$, $V_2=\pm e_1$ or $V_2=be_2+ce_3$ for any $b,c\in \R$, $b^2+c^2=1$ and $\displaystyle \phi(t)=c_1\left(\varepsilon \mu_1 t+\be\right)+\frac{c_2}{\left(\varepsilon \mu_1 t+\be\right)^3}, c_1,c_2\in \R$.
			
			\item 
			$\displaystyle G=SU(3),SO(3): \lambda_1=\lambda_2=\lambda_3>0, f(t)=\varepsilon\frac{\lambda_1}{\sqrt{2}}t+\be \text{ }(\varepsilon=\pm 1), \be\in\R, I=]-\frac{\be\sqrt{2}}{\lambda_1},+\infty[$ for $\varepsilon=1$, $\displaystyle I=]-\infty,\frac{\be\sqrt{2}}{\lambda_1}[$ for $\varepsilon=-1, V_2=ae_1+be_2+ce_3$, $a^2+b^2+c^2=1$  and \\ $\displaystyle \phi(t)=c_1(\varepsilon\frac{\lambda_1}{\sqrt{2}}t+\be)+\frac{c_2}{(\varepsilon \frac{\lambda_1}{\sqrt{2}}t+\be)^3}, c_1,c_2\in \R$.
			
					\end{enumerate}
	\end{prop}
	
	\begin{prop}[Non-unimodular groups cases]
		Let $G$ be a connected three-dimensional Riemannian non-unimodular  Lie group  equipped with an invariant metric and $V=V_1+V_2$  a harmonic vector field on the warped product $I\times_fG, I\subset \R$, $\displaystyle f: \R\longrightarrow ]0,+\infty[$ with $V_1=\phi(t)\partial_t$ and $V_2=ae_1+be_2+ce_3$ a unit left-invariant vector field on $G$. We have $$\displaystyle S(V_2)=\big[-\al^3(a^2+b^2)-\de^3(a^2+c^2)+\be(\al^2-\de^2)bc\big]e_1+a\big[\al^2\be c-\al\de^2 b+\be^2(\al-\de)b\big]e_2+a\big[-\be\de^2b-\al^2\de c-\be^2(\al-\de)c\big]e_3$$ and 
		$$\displaystyle \nabla_{V_2}{V_2}=(b^2\al+c^2\de)e_1+(-ac\be-ab\al)e_2+(ab\be-ac\de)e_3,$$ so $V$ is a harmonic map on $I\times_f G$  if and only if
		\begin{enumerate}
			
			\item $\al=\de>0$,$V_2=a e_1$, $a=\pm1$ and 
			\begin{eqnarray*}
				\begin{cases}\displaystyle
					f(t)f''(t)+2f'(t)=2\al^2,\\f''(t)\bigg(3f'(t)\phi^2(t)+f^2(t)f'(t)-2a\al f(t)\phi(t)\bigg)=0,\\ 4f'(t)^3a\phi(t) +\phi'(t)f(t)(a\al^2-2f'(t)^2)=2\al^3f(t)+2\al f(t)f'(t)^2 ,\\ 
					\phi''(t)+3\frac{f'(t)}{f(t)}\phi'(t)-3\left(\frac{f'(t)}{f(t)}\right)^2\phi(t)=-4a\al\frac{f'(t)}{f(t)}, t\in I .\end{cases}
			\end{eqnarray*}

			\item $\al>\de=0$,$V_2=a e_1$, $a=\pm1$ and 
			\begin{eqnarray*}
				\begin{cases}\displaystyle
					f(t)f''(t)+2f'(t)=2\al^2,\\f''(t)\bigg(3f'(t)\phi^2(t)+f^2(t)f'(t)-2a\al f(t)\phi(t)\bigg)=0,\\ -\al^3-\al f'(t)^2+4f''(t)f(t)\phi'(t)- a\phi(t)f'(t)f''(t)+\dfrac{\phi(t)f'(t)}{f(t)}(\al^2+2f'(t)^2)=0,\\ 
					\phi''(t)+3\frac{f'(t)}{f(t)}\phi'(t)-3\left(\frac{f'(t)}{f(t)}\right)^2\phi(t)=-2a\al\frac{f'(t)}{f(t)} , t\in I.\end{cases}
			\end{eqnarray*} 	
			\item $\al>\de>0$ or $\al>0>\de$ ,$V_2=a e_1+be_2+ce_3$, $a\ne 0$ and 
			\begin{eqnarray*}
				\begin{cases}\displaystyle	b(\be^2-\de^2)=\be(\al+\de)c,\\c(\be^2-\al^2)=-\be(\al+\de)b,\\
					f(t)f''(t)+2f'(t)=\al^2+\de^2,\\f''(t)\bigg(3f'(t)\phi^2(t)+f^2(t)f'(t)-a(\al+\de) f(t)\phi(t)\bigg)=0,\\ -\al^3(a^2+b^2)-\de^3(a^2+c^2)+\be bc(\al^2-\de^2)-af'(t)f''(t)\phi(t)+af(t)f''(t)\phi'(t)\\-a^2f'(t)^2(\al+\de)-f'(t)^2(b^2\al+c^2\de)+\dfrac{\phi(t)f'(t)}{f(t)}a(2f'(t)^2+\al^2+\de^2)=0,\\
					
					a(\al^2c \be-\al \de^2+\al\be^2-\de\be^2b)-bf'(t)\phi(t)f''(t)+f(t)\phi'(t)f''(t)b-abf'(t)^2(\al+\de)\\+af'(t)^2(c\be+b\al)+\dfrac{\phi(t)f'(t)}{f(t)}(2bf'(t)^2+b\al^2+b\al\de-c\al\be+c\be\de)=0,	\\

					a(-\de^2b\be-\al^2 \de c-\al\be^2c+\de\be^2c)-cf'(t)\phi(t)f''(t)+f(t)\phi'(t)f''(t)c-acf'(t)^2(\al+\de)\\-af'(t)^2(b\be-c\de)+\dfrac{\phi(t)f'(t)}{f(t)}(2cf'(t)^2+c\de^2+c\al\de-b\al\be+b\be\de)=0,
					\\ 
					\phi''(t)+3\frac{f'(t)}{f(t)}\phi'(t)-3\left(\frac{f'(t)}{f(t)}\right)^2\phi(t)=-2a(\al+\de)\frac{f'(t)}{f(t)}, t\in I.
				\end{cases}
			\end{eqnarray*} 	
		\end{enumerate}
	\end{prop}
	Using the examples of Propositions \ref{Eq: P1}, \ref{Eq: P2} and \ref{Eq: P3} for a non left-invariant vector fields on  Lie groups $G$, we obtain.
	\begin{prop}
		A harmonic vector field $V=\phi(t)\partial_t+a(x,y,z)\frac{\partial}{\partial x}$ on $I\times_f(\R\oplus\R\oplus\R)$ is a harmonic map if and only if $V_2$ is constant on $G$, $\phi$ is an affine function, and $f$ is a positive function on $I$ and $I=\R$.
		
	\end{prop}

	\begin{prop}
		A harmonic vector field $V=\phi(t)\partial_t+a(x,y,z)e_1$ on $I\times_f\mathbb{H}^3$ cannot be a harmonic map.
	\end{prop}
	
	\begin{prop}
		A harmonic vector field $V=\phi(t)\partial_t+b(y)e_2+c(y)e_3$ on  $I\times_f(\R\times H^2)$ cannot be a harmonic map.
	\end{prop}

\end{document}